\documentclass[onefignum,onetabnum]{siamart251216}

\usepackage{amsmath,amssymb,amsfonts,euscript,graphicx,hyperref,cleveref,mathtools,subcaption}
\usepackage{algorithm}
\usepackage{algpseudocode}
\usepackage[normalem]{ulem}
\usepackage{array}

\usepackage{placeins}

\def\bracket#1#2{\langle#1, #2 \rangle}

\newcommand{\tr}{\operatorname{tr}}
\newcommand{\trace}{\operatorname{tr}}
 
 \renewcommand{\ker}{\operatorname{Ker}}
 \renewcommand{\div}{\operatorname{div}}

\newcommand{\N}{\mathbb{N}}

\newcommand{\Nel}{{N_{\rm el}}}

\renewcommand{\H}{\mathcal{H}} % Hilbert space
\newcommand{\HS}[1]{\|#1\|_{\rm HS}} % Hilbert-Schmidt norm
\newcommand{\BH}{\mathcal{B}(\mathcal{H})} % Bounded operators on H

\newcommand{\Mn}{\mathcal{M}_n(\mathbb{C})} % Set of matrices size n
\newcommand{\D}{{\mathcal D}} % Set of density operators
\newcommand{\Dp}{{\mathcal D}^+} % Set of density operators
\newcommand{\Dn}{{\mathcal D}_n} % Set of density matrices of size n
\newcommand{\Dnp}{{\mathcal D}^+_n} % Set of POSITIVE density matrices of size n
\newcommand{\Herm}{\mathcal{H}_n(\mathbb{C})} % Set of Hermitian matrices size n
\newcommand{\Skewmat}{\mathcal{S}_n(\mathbb{C})} % Set of skew-hermitian matrices size n

\newcommand{\Barr}{\mathrm{Barr}} % Barrier function

\theoremstyle{plain}

\title{An algorithm for dynamical quantum optimal transport with applications to quantum chemistry}

\headers{Quantum transport applied to quantum chemistry}{G. Dusson, V.Ehrlacher, and E. Obermeyer}

\author{Genevi\`eve Dusson\thanks{Université Marie et Louis Pasteur, CNRS, LmB (UMR 6623), Besançon, France (\email{genevieve.dusson@math.cnrs.fr},  \email{etienne.obermeyer@math.cnrs.fr.})}
\and Virginie Ehrlacher\thanks{CERMICS, Institut Polytechnique de Paris, Ecole des Ponts \& Inria Paris, France, (\email{virginie.ehrlacher@enpc.fr}}).
\and \'Etienne Obermeyer\footnotemark[2]}

\date{\today}

\begin{document}

\maketitle

\begin{abstract}
Quantum optimal transport (QOT) is a rapidly developing field. Among the many formulations of this adaptation of classical optimal transport (OT) to spaces of density matrices, we numerically study a family of distances based on a dynamical formulation inspired by the Benamou--Brenier OT formulation.
We introduce an interior-point regularized method to compute geodesics between positive semidefinite matrices and visualize the results in terms of integral kernels and densities, inspired by quantum chemistry applications.
We show that dynamical QOT may provide a good approximation to certain problems in quantum chemistry with appropriate parameter tuning.
We also study the numerical properties of the distances at hand, and the convergence of the computed geodesics as the size of the matrices increases.
\end{abstract}



\section{Introduction} 

Optimal transport is a theory that defines a family of distances between probability measures called Wasserstein distances.
Introduced in 1781 by Monge~\cite{monge81}, it has since been largely developed, with major advances since the 1990s.
One of the most interesting features of Wasserstein spaces is the  existence of geodesics.
The Wasserstein metric can indeed be written in a Riemannian metric formalism, measuring the shortest path length between two measures. 
This fruitful point of view is due to~\cite{otto01}, developing a formula first introduced by Benamou and Brenier~\cite{benamou00} in 2000.
The geodesics which naturally correspond to barycenters between two measures with varying parameters have been generalized to  barycenters between more than two measures~\cite{agueh11}
and used in various fields, such as imaging science and machine learning~\cite{rabin12,Cuturi2013-fc, Bonneel2016-wd,Bigot2018-rk,dognin19}.
The resounding success of optimal transport theory has motivated the  development of  several extensions, such as unbalanced optimal transport~\cite{sejourne23}, Gromov--Wasserstein distance~\cite{Memoli2011-dl}, or modified Wasserstein distance on Gaussian mixtures~\cite{Delon2020-wk}.
In this article we are interested in Quantum Optimal Transport (QOT) theory, where the objects of interest are quantum probabilities~\cite{meyer95}.
The latter is a domain of functional analysis concerned with non-commutative structures that contains many extensions of traditional results in probability theory with general traces replacing integral calculus, operators or matrices replacing functions, density matrices replacing probability measures.
Many attempts at generalizing optimal transport theory in this direction have been made, most of them being based on various adaptations of the Kantorovich  formulation, which has the most convenient structure among the equivalent formulations of optimal transport as it is linear in the unknown~\cite{cole23,caglioti21,duvenhage22}. 
For a comprehensive overview and comparison of these various theories, we recommend~\cite{beatty25}.
In this work, we are interested in a formulation of QOT that would  allow the computation of barycenters between two or more density matrices. 
Because it is not clear whether Kantorovich-based formulations  give rise to a geodesic metric space, we choose to investigate the dynamical formulation introduced by Carlen and Maas~\cite{carlen20}.
They define a Riemannian metric on the set of density matrices  which naturally comes with geodesics.
Note that different properties of this metric have been investigated in a series of works by Chen et al.~\cite{chen16,chen17,chen17a}. 
For example in~\cite{chen17a}, a dual formulation for the problem that defines the metric is provided, and 
in~\cite{chen16} the distance is adapted to matrix-valued measures  by mixing classical and quantum optimal transport. 

The goal of  this work is to numerically investigate some properties of this dynamical QOT distance.
Computational classical optimal transport is now a well-developed field~\cite{peyre20}, and several 
numerical methods have been developed specifically for the dynamical formulation of optimal transport, starting with an augmented Lagrangian algorithm in the seminal paper by Benamou and Brenier~\cite{benamou00}. This algorithm is 
a particular case of proximal splitting methods detailed e.g. in~\cite{papadakis14} for the Benamou--Brenier problem, which
 are very well-adapted to this problem, as they can automatically handle vanishing densities.
Another possibility for the computation of dynamical optimal transport is given by the standard sequential quadratic programming (SQP) method, a smooth optimization method that resembles the Newton method for constrained optimization. 
It has been applied to dynamical optimal transport in~\cite{haber15} and 
in~\cite{natale21}, the latter using  an interior point method to  handle vanishing densities.

\begin{figure}[t]
    \centering
    \subcaptionbox{$t=0$}{\includegraphics[width=0.19\textwidth]{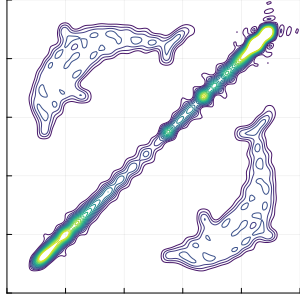}}
    \subcaptionbox{$t=0.25$}{\includegraphics[width=0.19\textwidth]{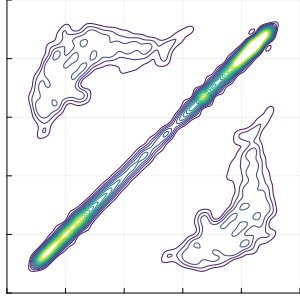}}
    \subcaptionbox{$t=0.5$}{\includegraphics[width=0.19\textwidth]{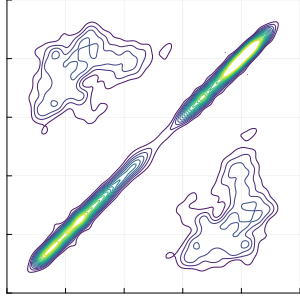}}
    \subcaptionbox{$t=0.75$}{\includegraphics[width=0.19\textwidth]{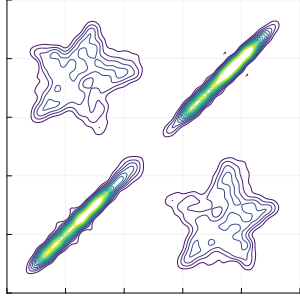}}
    \subcaptionbox{$t=1$}{\includegraphics[width=0.19\textwidth]{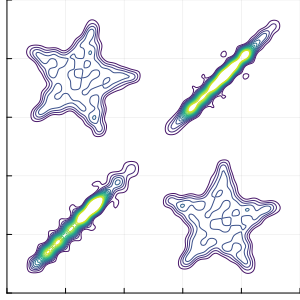}}
    \caption{A QOT geodesic between two density matrices (endpoints), represented by their kernels.}
    \label{fig:geod}
\end{figure}

\begin{figure}
    \centering
    % our geod
    \includegraphics[width=0.19\textwidth]{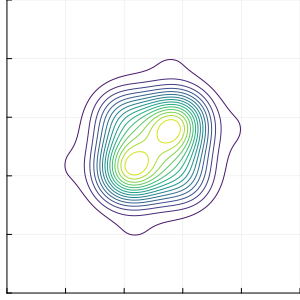}
    \includegraphics[width=0.19\textwidth]{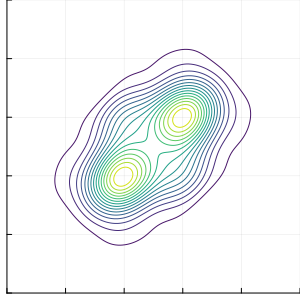}
    \includegraphics[width=0.19\textwidth]{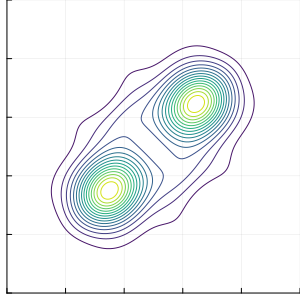}
    \includegraphics[width=0.19\textwidth]{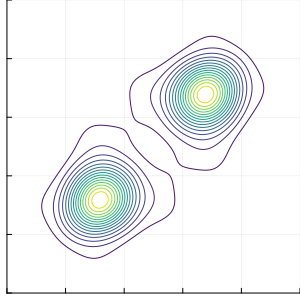}
    \includegraphics[width=0.19\textwidth]{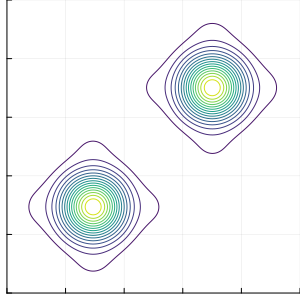}
    
    % chem curve
    \subcaptionbox{$t=0$}{\includegraphics[width=0.19\textwidth]{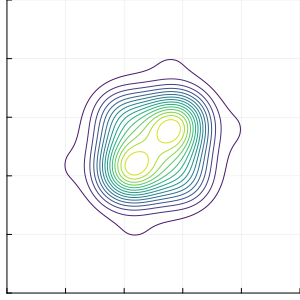}}
    \subcaptionbox{$t=0.25$}{\includegraphics[width=0.19\textwidth]{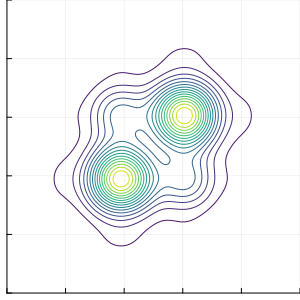}}
    \subcaptionbox{$t=0.5$}{\includegraphics[width=0.19\textwidth]{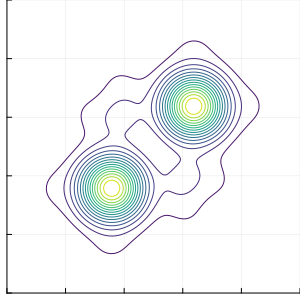}}
    \subcaptionbox{$t=0.75$}{\includegraphics[width=0.19\textwidth]{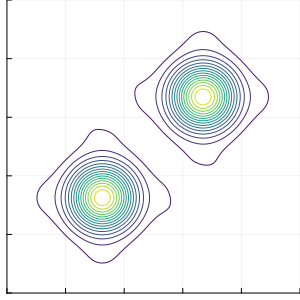}}
    \subcaptionbox{$t=1$}{\includegraphics[width=0.19\textwidth]{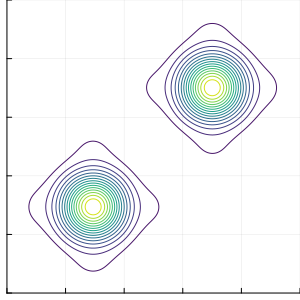}}
    \caption{Curves of density matrices represented by their kernels. Top: a  QOT geodesic computed with our algorithm (\Cref{sec:chem}). Bottom: electronic densities obtained by solving the electronic Schrödinger equation for varying nuclear positions (not a geodesic). Computations are made with SchrodingerFE.jl \cite{quan23}.}
    \label{fig:chem_geod}
\end{figure}

Numerical methods for dynamical QOT are so far limited, except for~\cite{chen17}, where Chen et al. have proposed a numerical algorithm for the practical computation of a matrix-valued Carlen--Maas distance, based on~\cite{haber15}. However in this work the distance was presented as a sum between a classical OT distance and QOT distance with a contribution of $99 \%$ classical and $1\%$ quantum distances in the presented numerical results.
In this work, we propose an algorithm for computing the standard Carlen--Maas QOT distance, by first 
adapting the algorithm presented in~\cite{chen17}
 for single density matrices (as opposed to matrix-valued measures), effectively computing for the first time, at least up to our knowledge, a fully quantum distance between density matrices, and their quantum interpolations. 
We also modify the algorithm to reach convergence with positive semidefinite matrices, by adding, in the spirit of~\cite{natale21}, a log-barrier term and a regularization strategy.
The last point notably allows us to carry out the computations for larger matrices, since large matrices with trace one often have small eigenvalues, as well as density matrices from quantum chemistry which are 
inherently of low rank. 
This also allows us to study different geodesics coming from different choices of the Carlen--Maas metric, and compare them to given density matrices from quantum chemistry. 
To give examples, we display
in~\cref{fig:geod} and~\cref{fig:chem_geod} (top row)  two instances of QOT geodesics between chosen density matrices (at $t=0$ and $t=1$). In \cref{fig:chem_geod} the density matrices provided on the bottom row are solutions to the electronic Schr\"odinger equation with varying nuclei positions, and we observe that the geodesics seem close to the one given by dynamical QOT (see \Cref{sec:chem}
for a more detailed comparison).

The main contributions of this article are the following. First, we adapt the
finite-volume SQP strategy developed in previous works on dynamical optimal
transport and matrix-valued transport to the computation of the standard
Carlen--Maas distance between single density matrices, thereby providing, to
the best of our knowledge, one of the first numerical methods for computing
fully quantum dynamical transport geodesics in this setting. Second, we make
this approach robust for positive semidefinite and low-rank density matrices by
combining the SQP scheme with a log-barrier interior-point method and a
trace-preserving endpoint regularization, which is essential for the quantum
chemistry examples considered here. Third, we investigate the numerical
behavior of the resulting geodesics with respect to the matrix dimension and
the time discretization, and we propose kernel- and density-based
visualizations that make it possible to compare dynamical QOT geodesics with
classical optimal transport interpolations. Finally, we apply the method to
families of density matrices arising from parametrized quantum chemistry
calculations, and we show that, after a suitable choice of the derivations
defining the Carlen--Maas metric, dynamical QOT geodesics can provide accurate
interpolations of such physically motivated curves of density matrices.

The outline of the article is as follows.
We begin in \Cref{sec:motivation} by 
recalling the definition of density matrices and how they arise in quantum chemistry, as well as the definition and convex formulation of dynamical QOT.
We then adapt~\cite{chen17} for the single matrix Carlen--Maas distance, recalling their discretization in \Cref{sec:discretization}.
In \Cref{sec:algo}, we recall the SQP algorithm used to solve the problem, and we discuss the adjunction of a barrier method and regularization to improve the performance of the algorithm for positive semidefinite matrices. 
We then discuss the practical implementation of the algorithm.
In \Cref{sec:numerics}, we perform a series of numerical experiments,  exploring the performance of the algorithm and displaying the results with kernels and density functions. 
This allows to visualize the properties of this distance for larger matrices.
We also compare the qualitative properties of dynamical QOT distances with the classical 2-Wasserstein distance.
\Cref{sec:chem} is dedicated to potential applications for quantum chemistry calculations where  density matrices naturally arise.
We visualize the results in this setting and discuss the possibility of constructing a reduced-order model based on QOT.
 Finally we provide concluding remarks in \Cref{sec:concl}.

\section{Setting: Density operators, density matrices, and Carlen--Maas Riemannian QOT distance} \label{sec:motivation}

In this section we present the framework of density operators and how they are approximated by density matrices. We then introduce the motivating quantum-chemistry problems, as well as the dynamical QOT formulation of Carlen and Maas~\cite{carlen20}. 

We denote by $\overline{x}$ the complex conjugate of $x$, and by $*$ the adjoint of either an operator or a matrix. In the context of complex Hilbert spaces, we use the convention that the inner product is linear in the second argument and antilinear in the first.

\subsection{Density operators}

We first define density operators and provide a few properties.
Let $\H$ be a separable Hilbert space, and let $\BH$ denote the space of bounded operators on $\H$.

\begin{definition}
An operator $T\in \BH$ is called a density operator on $\H$, if
\begin{enumerate}
    \item $T$ is self-adjoint and positive,
    \item $T$ is trace-class,
    \item $\operatorname{tr}(T)=1$.
\end{enumerate}
The set of density operators is denoted by $\D$. We also define the set of positive definite density operators,
$$\Dp = \{ T \in \BH, \ T^* = T, \ T \succ0, \ \tr(T)=1\}.$$
\end{definition}
Trace-class operators, as well as the few notions reminded here, are studied in \cite[VI.6]{reed80}.
We focus in this article on the Hilbert space $L^2([0,1])$ for simplicity.
\begin{proposition}
    Any density operator $T \in \D$ is a Hilbert--Schmidt (HS) operator. Thus, it is a kernel operator: there exists a function $\gamma \in L^2([0,1] \times [0,1])$, called the (integral) kernel of $T$, such that
    \begin{equation} \label{eq:HS}
    T = T_\gamma : \phi\in L^2([0,1]) \mapsto \left( y \mapsto \int_{[0,1]} \gamma(x,y) \phi(x)  \mathrm{d}x \right) \in L^2([0,1]).    
    \end{equation}
\end{proposition}
A natural norm on density operators is the Hilbert--Schmidt norm, defined by 
\[ 
\HS{T}:= \sqrt{\operatorname{tr}(T^*T)}.
\]
We also have $\HS{T_\gamma} = \|\gamma\|_{L^2([0,1]\times[0,1])}$.

\subsection{Density matrices}\label{subsec:density_matrices}

To approximate density operators, we use finite density matrices defined for a given dimension $n\in \N$ as 
\[
    \Dn = \{ \rho \in \Mn, \ \rho^* = \rho, \ \rho \succeq 0, \ \tr(\rho)=1\},
\]
where $\Mn$ is the space of $n \times n$ complex matrices. 
This approximation is justified by the fact that finite-rank operators are dense in $\D$ for (e.g.) the Hilbert--Schmidt norm. 
We approximate density operators $T_\gamma$ by density matrices denoted by $\rho$ through a Fourier  discretization of the kernel $\gamma$. However, any other discretization basis can be used.
More precisely, let 
\[(e_k : x \in [0,1] \mapsto \exp(2\mathrm i\pi k x))_{k\in \mathbb Z}\]
be the standard family of Fourier modes on $[0,1]$.
Let $\gamma \in L^2([0,1]\times [0,1])$ be a kernel corresponding to a density operator $T_\gamma$, and $K\in \N$ such that $2K+1 =n$. 
The approximate kernel 
\[
    \gamma_n \in {\rm Span}\left( \{ \overline{e_k} \otimes e_l, \; k,l \in \{ -K,\ldots, K\} \} \right) 
\]
is defined by 
\begin{equation}
\label{eq:gam_rho}
        \gamma_n : (x,y) \mapsto \sum_{k,l = -K}^K \rho^n_{k,l} \; \overline{e_k(x)} {e_l(y)},
\end{equation}
for $\rho^n\in \Dn$ defined by
\[
   \forall k,l \in \{ -K,\ldots, K\}, \
   \rho^n_{k,l}  = c \; \langle\overline{e_k}\otimes e_l ,\gamma \rangle
   = c \int_{[0,1]^2}  \hspace{-.2cm} \gamma(x,y){e_k(x)}\overline{e_l(y)} \ \mathrm dx \mathrm dy,
\]
with $c > 0$ ensuring that $\tr(\rho^n)=1$. 
Using this we have $T_{\gamma_n}
\xrightarrow[n\to+ \infty]{}
T_{\gamma}$ in Hilbert--Schmidt norm~\cite[VI.23]{reed80}. Note that for all $K\in \N$, using $n=2K+1$
\[
{\|\gamma_n\|_{L^2} = \|\rho^n\|_2},
\]
where $\|\rho^n\|_2$ denotes the Frobenius norm of matrices.

\subsection{Motivation: Density matrices in quantum chemistry} \label{subsec:motivation}

In many quantum chemistry models, the central object is often not the full many-electron wavefunction but the electronic density matrix.
It arises naturally when one seeks a more compact representation of a quantum state: the full wavefunction solution to the electronic time-independent Schr\"odinger equation lives in a very high-dimensional space and contains far more information than needed to compute physical observables. 
The density matrix reduces this complexity by projecting the many-body state onto the one-particle space; it encodes the electronic distribution, essential correlations, and allows energies and expectation values to be written as matrix traces.
More precisely, given a wavefunction 
$\psi \in L^2(\mathbb [0,1]^{\Nel})$ where $\Nel$ corresponds to the number of electrons in the system, the associated kernel $\gamma \in L^2([0,1]\times [0,1])$ is defined by
\[
    \gamma(x,y) = \int_{[0,1]^{\Nel-1}} \overline{\psi(x,x_2,\dots,x_\Nel)}\psi(y,x_2,\dots,x_\Nel) \ \mathrm d x_2 \dots \mathrm dx_\Nel,
\]
and the associated density operator is the integral operator $T_\gamma$, defined as in equation~\cref{eq:HS}.
A typical aim in quantum chemistry is the efficient computation of kernels $\gamma_\theta$ for different values of a parameter $\theta$ corresponding e.g. to varying positions of nuclei. 
It is then necessary to carry out demanding  calculations for a large number of parameters~$\theta$.
A possible way to tackle this problem is to use interpolation techniques between two or more values of $\theta$ and corresponding $\gamma_\theta$.
These can typically be obtained using linear reduced-order modeling techniques~\cite{Barrault2004-ys,Hesthaven2016-xm}, but such methods do not work well in general  when transport effects are playing a large role.
In the latter case, if kernels $\gamma_\theta$ were probability densities, one could compute optimal transport interpolations between them as proposed by some of us in~\cite{Dalery2026-pl}.
Unfortunately they are not, since kernels can be complex-valued and integrate to different numbers.
They only satisfy
\begin{equation*}
    \gamma_\theta(x,y) = \overline{\gamma_\theta(y,x)}, \quad \int_{[0,1] }\gamma_\theta(x,x) \ \mathrm{d} x = 1. 
\end{equation*} 
However, kernel operators $T_{\gamma_\theta}$ are density operators, i.e.
the counterpart in quantum probability of classical probability measures.
Thus, to build interpolations we consider geodesics for a suitable distance between density matrices described below.

\subsection{Dynamical QOT Riemannian distance}\label{sec:defs}

In~\cite{carlen20}, Carlen and Maas have proposed a dynamical QOT distance between density matrices.
The construction relies on a differential structure~\cite[Definition 4.7]{carlen20} and a theory of two-operator functional calculus~\cite[§6]{carlen20}.
The authors leave many choices of parameters which all lead to the construction of a Riemannian distance on the set of density matrices.
This distance is defined on the finite-dimensional von Neumann algebra $\Mn$ of $n\times n$ complex matrices, with the usual trace $\tr$.
Note that the construction is only available in a finite-dimensional setting, although notation adapts to infinite dimension. 
Let $\Herm$ denote the set of Hermitian matrices.
We define $J$ derivations denoted by $\partial_j$ for $j \in \{1, \dots, J\}$ acting as commutators with matrices $L_j \in \Herm$ for $j \in \{1, \dots, J\}$
\[
    \forall A\in \Mn, \ \partial_j A \vcentcolon = L_j A - A L_j. 
\]
The gradient operator $ \nabla_L :\Mn \to \Mn^J $  is defined by
\begin{equation} \label{eq:nabla}
    \nabla_L A = (\partial_jA)_{j\in \{1,\dots,J\}} = (L_j A - A L_j)_{j\in \{1,\dots,J\}}.
\end{equation}
The adjoint operator $\nabla^*_L$, also denoted $- \div_L$, acts on vectors of $J$ matrices $\mathbf{B} \in \Mn^J,$ according to
\begin{equation*}
    \nabla^*_L(\mathbf B) = -\sum_{j=1}^J \partial_j B_j. 
\end{equation*}
In order to build a distance defined on the whole set of positive-definite density matrices
\[
    \Dnp := \{ \rho \in \Mn, \ \rho^* = \rho, \ \rho \succ 0, \ \tr(\rho)=1\},
\]
the assumption $\ker(\nabla_L) = \operatorname{Span}(I_n)$ is imposed~\cite[Remark 7.1]{carlen20}. 
A function $\theta$ must also be chosen to define the action of density matrices on matrices~\cite[§6]{carlen20}. 
We adopt the simplest choice $\theta(a,b) = \frac{a+b}{2}$, which defines an action of a density matrix $\rho$ on a matrix $A$ by $\frac{1}{2}(\rho A + A \rho)$. 
Altogether this constructs a distance on the set $\Dnp$, given by,
for $\rho_0,\rho_1 \in \Dnp$,
 \begin{align}
	\mathcal{W}(\rho_0,\rho_1)^2 = &\min_{\rho_t\in \Dnp, \ \mathbf v_t\in \Skewmat^J} \int_0^1
	\trace(\rho_t \mathbf v^*_t \mathbf v_t)  \ \mathrm dt,  \label{eq:continuous_problem_v} \\
	&\frac{d \rho_t}{d t}-\frac{1}{2} \nabla_L^* (\mathbf v_t\rho_t+\rho_t \mathbf v_t)=0, \label{eq:continuous_problem_v_2} \\
    & {\rho_t}_{|t=0} = \rho_{0}, \quad {\rho_t}_{|t=1} = \rho_{1}. \label{eq:continuous_problem_v_3} 
\end{align}
The minimum is taken over smooth curves $(\rho_t)_{t\in [0,1]} \subset \Dnp$ and $(\mathbf v_t)_{t \in[0,1]} \subset \Skewmat^J$ satisfying the constraint conditions, where 
\[
    \Skewmat = \{A \in \Mn, \  A^* = -A \}
\]
denotes the set of skew-Hermitian matrices. 
Note that $\mathbf v_t$ is a column vector of matrices, and $\mathbf{v}_t^*$ is the row vector of adjoints, so that $\mathbf v_t^* \mathbf{v}_t = \sum_{j=1}^J ( \mathbf v^*_{t})_j (\mathbf v_{t})_j$. Since components of $\mathbf{v}_t$ are skew-Hermitian, the objective function is equal to
\begin{equation}\label{eq:continuous_objective_v}    
  \int_0^1
	\trace(\rho_t \mathbf v^*_t \mathbf v_t)  \ \mathrm dt =  -\int_0^1 \sum_{j=1}^J \tr(\rho_t (\mathbf v_t)_j^2) \ \mathrm dt. 
\end{equation}
Note that the dynamical QOT metric defined on $\Dnp$ can in fact be extended continuously to $\Dn$
as shown in~\cite[Proposition 9.2]{carlen20}.
Moreover a convex formulation is available for this problem. The variable \( \mathbf u_t = \mathbf v_t \rho_t \) is introduced, and $\overline{\mathbf{u}_t}$ denotes the column vector of adjoints $\overline{\mathbf{u}_t} = (\mathbf{u}_t^*)^T =  ((\mathbf{u}^*_t)_j)_{j\in \{1,\dots, J\}}$,
from which we easily obtain $\overline{\mathbf u_t} = -\rho_t \mathbf{v}_t$.
Using the cyclicity of the trace, \Crefrange{eq:continuous_problem_v}{eq:continuous_problem_v_3} become
\begin{align}
\mathcal{W}(\rho_0,\rho_1)^2 = & \min_{\rho_t\in \Dnp, \ \mathbf u_t\in \Mn^J} \int_0^1  \mathrm{tr}\left( \mathbf u_t^* \mathbf u_t \rho_t^{-1} \right) \, \mathrm{d} t \label{eq:continuous_problem_u} \\
& \frac{d \rho_t}{dt} -  \frac{1}{2} \nabla_L^* \left( \mathbf u_t - \overline{\mathbf u_t} \right) =0  \label{eq:continuity_equation}\\ 
& {\rho_t}_{|t=0} = \rho_{0}, \quad {\rho_t}_{|t=1} = \rho_{1}. \label{eq:continuous_problem_u_3}
\end{align}
Note that in~\cite{carlen20} the problem is posed with another set of variables, namely $(\rho_t,A_t)\in \Dn^+ \times \Mn$, with $\tr(A_t) = 0$, $A_t \in \Herm$, and $\mathbf v_t = \nabla_L A_t$.
In~\cite[3.2.1]{chen16}, it is proven that the optimal $\mathbf v_t$ in \crefrange{eq:continuous_problem_v}{eq:continuous_problem_v_3} 
is a gradient, namely
\[
\mathbf v_t = \nabla_L A^{\mathrm{cont}}_t,
\]
where $-A^{\mathrm{cont}}_t\in \Herm$  is a Lagrange multiplier associated to the continuity equation constraint~\eqref{eq:continuous_problem_v_2}. This multiplier can be chosen traceless, which shows that~\crefrange{eq:continuous_problem_v}{eq:continuous_problem_v_3} is equivalent to the problem with variables $(\rho_t,A_t)$~\cite{carlen20}.
Besides, in~\cite[Theorem 16]{chen17a}, it is shown that the optimal $\mathbf u_t$ in~\crefrange{eq:continuous_problem_u}{eq:continuous_problem_u_3} writes $\mathbf u_t = (\nabla_L A_t) \rho_t$, for some traceless $A_t\in\Herm$. The proof uses the fact that the representation $\mathbf u = (\nabla_L A) \rho$ holds for pairs $(\rho,\mathbf u)$ that minimize an  analogue of the kinetic energy given by the objective function \cref{eq:continuous_objective_v}.
The combination of the arguments shows that~\crefrange{eq:continuous_problem_v}{eq:continuous_problem_v_3} is indeed equivalent to~\crefrange{eq:continuous_problem_u}{eq:continuous_problem_u_3}.

\paragraph{Gradient choice}\label{sec:gradient1}
In order to finalize the distance definition we need to choose derivations and therefore corresponding operators $L_j \in \Herm$.
We tested a family of two-matrix sets parametrized by two values $\alpha \in \mathbb R_+^*,\  \beta\in \mathbb{R} \cup \{-\infty\}$,
\begin{align}\label{eq:derivations1}
     L_1(\alpha) =& \begin{pmatrix}
        -(\frac{n-1}{2})^\alpha & & (0) \\
        & \ddots & \\
        (0) &  & (\frac{n-1}{2})^\alpha
    \end{pmatrix}, \\
    L_2(\beta) =&
\begin{pmatrix}
0 & 1 & 2^\beta & \cdots & (n-1)^\beta \\
1 & 0 & 1 & \cdots & (n-2)^\beta \\
2^\beta & 1 & 0 & \cdots & (n-3)^\beta \\
\vdots & \vdots & \vdots & \ddots & \vdots \\
(n-1)^\beta & (n-2)^\beta & (n-3)^\beta & \cdots & 0
\end{pmatrix}.\label{eq:derivations2}
\end{align}
The first matrix $L_1(\alpha)$ is diagonal, with diagonal equal to $(\operatorname{sgn}(k) |k|^\alpha)_{k=-\frac{n-1}{2}}^{\frac{n-1}{2}}$,
while $L_2(\beta)$ is a Toeplitz matrix. Any set of two matrices $L = \{L_1(\alpha),L_2(\beta)\}$ satisfies the condition $\ker(\nabla_L) = \operatorname{Span}(I_n)$.
Unless otherwise stated, all numerical experiments up to~\Cref{sec:chem} are performed with $\alpha = 1$ and $\beta= -\infty$.

\section{Discretization}\label{sec:discretization}
In this section, we adapt the discretization proposed in~\cite{chen17} for  single matrices (as opposed to matrix-valued measures), which is based on a finite-volume discretization. 

\subsection{Discretization of the constraints}
The interval $[0,1]$ appearing in ~\crefrange{eq:continuous_problem_v}{eq:continuous_problem_v_3} and called, for simplicity, the time interval is discretized into $P$ intervals of equal size. Denote the timestep \( h = \frac{1}{P} \).
The curve of matrices $\rho$ is discretized at the boundary of the intervals, and the endpoints are excluded. 
The variable $\mathbf{u}$ is discretized at the centers of the intervals, following a staggered grid framework, which is best suited for this particular problem~\cite{papadakis14}.

The discretized variables are still denoted $\rho$ and $\mathbf{u}$ and read
\begin{equation}
    \rho = (\rho_{p+\frac12})_{p\in \{1,\dots, P-1\}}, \quad \mathbf{u} = (\mathbf u_p)_{ p \in \{1,\dots,P\}}.
\end{equation}
The discretized version of the continuity equation~\cref{eq:continuity_equation} is
\begin{equation}
D_1 \rho + D_2 \mathbf{u} = b,
\end{equation}
where the operators 
\[
    D_1 : (\Herm)^{P-1} \to (\Herm)^{P}\text{ and } D_2 :((\Mn)^{J})^{P} \to (\Herm)^{P}
\] 
are defined as
\[
(D_1 \rho)_p = 
\begin{cases}
\displaystyle \phantom{-}\frac{1}{h} \rho_{\frac{3}{2}} ,& p = 1 \\
\displaystyle\frac{\rho_{p+\frac{1}{2}} - \rho_{p-\frac{1}{2}}}{h} ,& 2 \leq p \leq P - 1, \\
\displaystyle - \frac{1}{h} \rho_{P - \frac{1}{2}} ,& p = P
\end{cases}
\quad
(D_2 \mathbf{u})_p = -\frac{1}{2} \nabla_L^*(\mathbf u_p - \overline{{\mathbf{u}}_p}).
\]
The right-hand side vector \( b \in (\Herm)^{P} \) encodes boundary information, namely
\[
b_p = 
\begin{cases}
\phantom{-}\frac{1}{h} \rho_0 ,&  p = 1 \\
-\frac{1}{h} \rho_1 ,& p = P \\
\phantom{-} 0 & \text{otherwise.}
\end{cases}
\]
\subsection{Discretization of the objective function}
The objective function in problem~\cref{eq:continuous_problem_u}
$
\displaystyle\int_0^1  \mathrm{tr}\left( \mathbf u_t^*  \mathbf u_t  \rho_t^{-1} \right) \, \mathrm{d} t
$
is discretized, following~\cite{chen17} and~\cite{haber15},   as
\begin{equation}\label{eq:discretized_obj_h}
 \sum_{p=1}^{P}  \text{tr}\left( \mathbf{u}_p^* \mathbf{u}_p  \left( \rho_{p - \frac{1}{2}}^{-1} + \rho_{p + \frac{1}{2}}^{-1} \right) \right)  \frac{h}{2}.
\end{equation}
This can be written using an averaging operator \( A : {(\Dn^+)}^{P-1} \hspace{-.25cm} \to (\Mn)^{P} \)  and a vector $a \in (\Mn)^{P}$ accounting for the boundary values \( \rho_0 \) and \( \rho_1 \), defined by
\[
(A(\rho^{-1}))_p = 
\begin{cases}
\frac{1}{2} \rho_{\frac{3}{2}}^{-1} ,& p = 1 \\
\frac{1}{2} \left( \rho_{p - \frac{1}{2}}^{-1} + \rho_{p + \frac{1}{2}}^{-1} \right) ,& 2 \leq p \leq P - 1 \\
\frac{1}{2} \rho_{P - \frac{1}{2}}^{-1} ,& p = P
\end{cases},
\ 
a_p = 
\begin{cases}
\frac{1}{2} \rho_0^{-1} ,& p = 1 \\
\frac{1}{2} \rho_1^{-1} ,& p = P \\
0 & \hspace{-.8cm} \text{otherwise.}
\end{cases}
\]
The objective function~\cref{eq:discretized_obj_h} becomes
\begin{equation*}
  h  \left\langle \mathbf{u}^* \circ \mathbf{u},\; A(\rho^{-1}) + a \right\rangle,
\end{equation*}
where $\circ$ is the block-wise product, that is, 
 \[
 \mathbf{u}^* \circ \mathbf{u}  = (\mathbf{u}^*_p  \mathbf{u}_p)_{p\in \{1,\dots,P\}} = 
 \left(\sum_{j=1}^J (\mathbf u_p)_j^* (\mathbf u_p)_j \right)_{p\in \{1,\dots,P\}},
 \] 
and
the scalar product on vector of matrices is defined as 
\[
\displaystyle\bracket{(A_p)_p}{(B_p)_p} = \sum_{p=1}^P \tr(A_p^*B_p).
\]
The discretized problem of \crefrange{eq:continuous_problem_u}{eq:continuous_problem_u_3} can therefore be expressed in the following way, where we discarded $h$ to obtain an equivalent formulation,
\begin{align} \label{eq:discretized_pb}
\min_{\rho \in (\Dnp)^{P-1}, \mathbf{u} \in (\Mn^J)^P} \quad &  \left\langle \mathbf{u}^* \circ \mathbf{u},\; A(\rho^{-1}) + a \right\rangle \\
\text{subject to } \quad & D_1 \rho + D_2 \mathbf{u} = b. \label{eq:discretized_pb2}
\end{align}

\section{Optimization algorithm}\label{sec:algo}
In this section, we present the resolution of the discretized problem~\crefrange{eq:discretized_pb}{eq:discretized_pb2}
using an SQP algorithm, following~\cite{chen17,haber15}, and improving it by proposing an interior point method as well as a regularization strategy.

\subsection{SQP formulation}
\label{subsec:SQP}
The SQP method amounts to using a Newton method on the gradient of the Lagrangian to solve the discretized problem. For an introduction, see for instance~\cite[part III]{bonnans10}. 
In this context, the Lagrangian function $\mathcal{L} :(\Dnp)^{P-1}\times (\Mn^J)^P \times (\Herm)^P  \to  \mathbb R$ is given by
\[ \mathcal{L}(
  \rho, \mathbf{u}, \lambda) =\left\langle \mathbf{u}^* \circ \mathbf{u},\; A(\rho^{-1}) + a \right\rangle + \left\langle \lambda,\; D_1 \rho + D_2 \mathbf{u} - b \right\rangle,
\]
and the Karush-Kuhn-Tucker (KKT) conditions are
\begin{align*}
\nabla_\rho \mathcal{L} &= D_1^* \lambda -   \rho^{-1} \circ A^*(\mathbf{u}^* \circ \mathbf{u}) \circ \rho^{-1} = 0, \\
\nabla_u \mathcal{L} &= D_2^* \lambda + 2  \mathbf{u} \circ (A(\rho^{-1}) + a) = 0, \\
\nabla_\lambda \mathcal{L} &= D_1 \rho + D_2 \mathbf{u} - b = 0.
\end{align*}
The Hessian of the objective function ${f}(\rho,\mathbf{u}) 
:= \left\langle \mathbf{u}^* \circ \mathbf{u},\; A(\rho^{-1}) + a \right\rangle$ 
is given by
\begin{equation}
    \label{eq:discrete_hessian}
\hat H = 
    \begin{pmatrix}
        \nabla_\mathbf{\rho} \nabla_\mathbf{\rho}  {f}(\rho, \mathbf{u}) & (\nabla_\mathbf{\rho} \nabla_\mathbf{u}  {f}(\rho, \mathbf{u}))^* \\
        \nabla_\mathbf{\rho} \nabla_\mathbf{u}  {f}(\rho, \mathbf{u}) & \nabla_\mathbf{u} \nabla_\mathbf{u}  {f}(\rho, \mathbf{u})
    \end{pmatrix},
\end{equation}
and a possible block diagonal approximation is
\begin{align}  
    &\hat H = \label{eq:discrete_hessian_blockdiag}
    \begin{pmatrix}
        \nabla_\mathbf{\rho} \nabla_\mathbf{\rho}  {f}(\rho, \mathbf{u}) & 0 \\
        0 & \nabla_\mathbf{u} \nabla_\mathbf{u}  {f}(\rho, \mathbf{u})
    \end{pmatrix}.
\end{align}
We provide the explicit computation of these matrices in Appendix~\ref{app:Hessian}.
At each SQP iteration, the following saddle point linear system is solved,
\begin{equation}\label{eq:saddle_point}
\begin{bmatrix}
\hat{H} & D^* \\
D & 0
\end{bmatrix}
\begin{bmatrix}
\delta w \\
\delta \lambda
\end{bmatrix}
= -
\begin{bmatrix}
\nabla_w \mathcal{L} \\
\nabla_\lambda \mathcal{L}
\end{bmatrix},
\end{equation}
with $\hat{H}$ given either by~\cref{eq:discrete_hessian} or \cref{eq:discrete_hessian_blockdiag}, \( w = (\rho, \mathbf{u}) \) is the vector of primal variables and \( D = (D_1, D_2) \) is the constraint operator. 
This is a symmetric linear system. Note that we implemented the spaces $\Herm$ and $\Mn$ as real vector spaces. Their dimensions are $\dim_{\mathbb{R}}(\Herm) = n^2$ and $\dim_{\mathbb{R}}(\Mn)= 2n^2$, so that the size of the linear system~\cref{eq:saddle_point} is
\[
   \underbrace{n^2(P-1)+2n^2PJ}_{\text{size of $w=(\rho,\mathbf u)$}}+\underbrace{n^2P}_{\text{size of $\lambda$}}  = 2Pn^2(J+1) - n^2 .
\]
The variables are then updated by
\[(w,\lambda) \xleftarrow{} (w,\lambda) + \eta(\delta w, \delta \lambda),\]
where $\eta \in (0,1]$ is determined by a line search. In the implementation we chose to use Bonnans et al.~\cite[§17, §18]{bonnans10}.
We summarize the SQP method in~\cref{alg:SQP}.
\begin{algorithm}[H]
\caption{SQP Method : \texttt{SQP(}$\mathcal{L},\hat{H},D,(\rho,\mathbf{u},\lambda),\tau$\texttt{)}}
\label{alg:SQP}
\begin{algorithmic}[1]
\Require
\Statex \begin{itemize}
    \item Lagrangian function $\mathcal{L}$,
    \item Hessian approximation $\hat{H}$, constraint operator $D$
    \item initial guesses $\rho \in (\Dnp)^{P-1}$, $\mathbf{u} \in (\Mn^J)^P$, $\lambda\in (\Herm)^P$, 
    \item tolerance $\tau>0$.
\end{itemize}
\State Compute the algorithm error $\mathrm{err} = \|\nabla \mathcal{L}(\rho,\mathbf{u},\lambda)\|_2$
\While{$\mathrm{err} > \tau$}
    \State Compute $ (\delta\rho, \delta \mathbf{u}, \delta \lambda)$ by solving~\cref{eq:saddle_point}
    \State Compute an appropriate step $\eta>0$ via line search on a merit function
    \State Update $(\rho,\mathbf{u},\lambda) = (\rho,\mathbf{u},\lambda) + \eta (\delta\rho, \delta \mathbf{u}, \delta \lambda) $
    \State Update the algorithm error $\mathrm{err} = \|\nabla \mathcal{L}(\rho,\mathbf{u},\lambda)\|_2$
\EndWhile
\end{algorithmic}
\end{algorithm}

In our implementation, we have used the following initial guesses for the SQP algorithm.
\begin{enumerate}
    \item The vector of density matrices $\rho \in (\Dnp)^{P-1}$ is initialized by the linear combinations $\rho_p = \frac{P+1-p}{P+1} \rho_0 + \frac{p}{P+1} \rho_1, \ p\in \{1,\dots,P-1\}$.
    In our experiments, this guess has proven to be close to the solution.
    \item The variable $\mathbf{u} \in (\Mn^J)^P$ is initialized as $\rho \mathbf v$, with $\rho$ defined above, and $\mathbf v  = \nabla_L \left(\frac{MM^*}{\tr(MM^*)}\right)$ the gradient of a matrix $M \in \Mn$ with entries sampled from a standard normal distribution.
    \item The dual variable $\lambda\in (\Herm)^P$ is initialized as a vector of identity matrices.
\end{enumerate}

\subsection{Interior point method and regularization}

A major difficulty arising in practice is that~\cref{alg:SQP} may not converge, or may converge slowly, if the density matrices are not positive definite or when their eigenvalues are too small, which is often the case in practice. Actually, without addressing this problem, many of the computations we tried would simply not converge.
The same issue also arises in the classical optimal transport setting with vanishing densities and can be solved with different strategies, e.g.,~\cite{natale21} using an interior-point method, and~\cite{papadakis14} using a first-order proximal splitting method, which converges unfortunately too slowly for the current purpose. 

\paragraph{Interior point method}
We therefore  adapt~\cite{natale21} and   propose an interior point method for this problem.
To do so, we modify the objective function and add the barrier function
 \begin{equation*}
     \Barr:\rho  \mapsto 
    \begin{cases} \displaystyle
      - \sum_{p=1}^{P-1} \ln(\det(\rho_{p + \frac12})) & \text{ if } \, \forall p \in \{1,\dots,P-1\}, \; \rho_{p + \frac12}  \succ 0,   \\
      +\infty & \text{ else.} 
    \end{cases}
\end{equation*}
 This is the classical log barrier function in positive semidefinite optimization (see e.g.~\cite[§6.4]{nesterov94}). 
 The time-discretized problem becomes
\begin{align}\label{eq:pb_barrier}
\min_{\rho \in (\Dnp)^{P-1},\  \mathbf{u} \in (\Mn^J)^{P}} \quad & \left\langle \mathbf{u}^* \circ \mathbf{u},\; A(\rho^{-1}) + a \right\rangle + \mu  \Barr(\rho) \\
\text{subject to } \quad & D_1 \rho + D_2 \mathbf{u} = b,
\label{eq:pb_barrier2}
\end{align}
with $\mu \in [0,+\infty).$ This is a linearly constrained convex problem. 
The Lagrangian   becomes
\begin{equation}
\label{eq:Lagrangianmu}
    \mathcal{L}_\mu:(\rho, \mathbf{u}, \lambda) \mapsto \left\langle \mathbf{u}^* \circ \mathbf{u},\; A(\rho^{-1}) + a \right\rangle +\mu \Barr(\rho) + \left\langle \lambda,\; D_1 \rho + D_2 \mathbf{u} - b \right\rangle.
\end{equation}
Noting that 
\[
    \nabla_\rho \Barr(\rho)=-\rho^{-1},
\]
where the inverse is to be understood componentwise, 
the KKT conditions become
\begin{align*}
\nabla_\rho \mathcal{L}_\mu &= D_1^* \lambda - \mu \rho^{-1} - \rho^{-1} \circ A^*(\mathbf{u}^* \circ \mathbf{u}) \circ \rho^{-1} = 0, \\
\nabla_u \mathcal{L}_\mu &= D_2^* \lambda + 2  \mathbf{u} \circ (A(\rho^{-1}) + a) = 0, \\
\nabla_\lambda \mathcal{L}_\mu &= D_1 \rho + D_2 \mathbf{u} - b = 0,
\end{align*}
with 
\begin{equation}
\label{eq:D1star}
    (D_1^* \lambda)_{p=1}^{P-1} = \left(\frac{\lambda_p - \lambda_{p+1}}{h}\right)_{p=1}^{P-1},\quad
(D_2^* \lambda)_{p=1}^{P} 
= -\frac{1}{2}( \nabla_L\lambda_p - \overline{\nabla_L \lambda_p})_{p=1}^P.
\end{equation}
We then solve problem~\crefrange{eq:pb_barrier}{eq:pb_barrier2} for different values of $\mu,$ as $\mu \to 0$, using intermediate results for warm starts.
Note that in standard (linear) semidefinite programming, it is often best to lower the value of $\mu $ at each iteration~\cite{nesterov94} of the Newton method in the context of the interior point method. However, in our computations, this strategy sometimes resulted in the divergence of the algorithm, so we chose to decrease $\mu$ when an intermediate tolerance was reached.

\paragraph{Regularization}

In practice the barrier function was not always sufficient to guarantee the convergence of the algorithm. We therefore combined it with a regularization method, which is another method one may use to tackle the poor behavior of the SQP algorithm when the density matrices are singular.
It consists of computing the QOT distance between $\rho_0^\varepsilon$ and $\rho_1^\varepsilon$ defined as
\begin{equation}
\label{eq:rhoeps}
    \rho_0^\varepsilon := (1+n\varepsilon)^{-1}(\rho_0 + \varepsilon I_n) \quad 
    \rho_1^\varepsilon := (1+n\varepsilon)^{-1}(\rho_1 + \varepsilon I_n),
\end{equation}
for varying values of $\varepsilon$ decreasing to $0$. 
The computations are expected to converge to the distance between $\rho_0$ and $\rho_1$, as shown in~\cite[§9]{carlen20}.

The resulting algorithm that combines all the modifications of the SQP algorithm is summarized in~\cref{alg:descent}. 
We use the function \texttt{SQP()} from~\cref{alg:SQP}.

\begin{algorithm}[H]
\caption{Barrier--Regularization Descent Scheme}
\label{alg:descent}
\begin{algorithmic}[1]
\Require Initial data: $\rho_0 \in \Dn, \rho_1 \in \Dn$
\Require Initial and target values for barrier $\mu_\text{init}$, $\mu_\text{end}$; regularization $\varepsilon_\text{init}$, $\varepsilon_\text{end}$ and tolerance: $\tau_\text{init}$, $\tau_\text{end}$.
\State Initialize barrier $\mu = \mu_\text{init}$, regularization $\varepsilon = \varepsilon_\text{init}$ and  tolerance $\tau = \tau_\text{init}$
\State Compute the trace-one regularized $\rho^\varepsilon_0$, $\rho_1^\varepsilon$ using~\cref{eq:rhoeps}
\State Compute the initial guesses $(\rho, \mathbf{u}, \lambda)$ as in~\Cref{subsec:SQP}
\While{$\mu \geq \mu_\text{end}$ \textbf{or} $\varepsilon > \varepsilon_\text{end}$}
    \State Build the Lagrangian $\mathcal{L}^\varepsilon_\mu$ using \cref{eq:Lagrangianmu} (depends on $\varepsilon$ through $\rho^\varepsilon_0, \rho^\varepsilon_1$ in $a,b$)
    \State Build the Hessian $\hat{H}^\varepsilon_\mu$ using~\cref{eq:discrete_hessian} or~\cref{eq:discrete_hessian_blockdiag}
     (depends on $\varepsilon$ through $\rho^\varepsilon_0, \rho^\varepsilon_1$)
    \State $(\rho,\mathbf{u},\lambda) \gets $ \texttt{SQP(}$\mathcal{L_\mu^\varepsilon},\hat{H}_\mu^\varepsilon,D,(\rho,\mathbf{u},\lambda),\tau$\texttt{)}
    \If{$\varepsilon > \varepsilon_\text{end}$}
        \State Decrease the regularization $\varepsilon$, e.g. $\varepsilon \gets \varepsilon/10$
        \State Compute the new trace-one regularized $\rho^\varepsilon_0,\rho^\varepsilon_1$
    \Else
        \State Decrease the barrier $\mu$, e.g. $\mu \gets \mu/10$
        \State Decrease the tolerance $\tau$, e.g. $\tau \gets \tau/2$ 
    \EndIf
\EndWhile
\State \textbf{Return} $(\rho,\mathbf{u},\lambda)$
\end{algorithmic}
\end{algorithm}

\subsection{Solving the linear system}

The bottleneck of the SQP method is to compute the solution to the linear system~\cref{eq:saddle_point}. 
Several articles are dedicated to the topic, see for instance~\cite{benzi05,benzi08}.
Two main methods can be used, either solving the linear system directly or iteratively, the latter typically relying on GMRES (Generalized minimal residual), or using a Schur complement method as in~\cite{chen17}, 
by first solving the reduced system
\begin{equation}\label{eq:schur}
    D \hat{H}^{-1} D^* \delta \lambda = \nabla_\lambda \mathcal{L} - D \hat{H}^{-1} \nabla_w \mathcal{L},
\end{equation}
and then compute the primal descent direction as
\[
\delta w = -\hat{H}^{-1} (D^* \delta \lambda + \nabla_w \mathcal{L}).
\]
In both cases, the exact Hessian~\cref{eq:discrete_hessian} or an approximation thereof such as~\cref{eq:discrete_hessian_blockdiag} can be used.

For the Schur complement a first difficulty is that the Schur matrix $ D \hat{H}^{-1} D^*$ is not positive-definite because $D$ is not surjective since 
\[
    \ker(D \hat{H}^{-1} D^*) = \ker(D^*),
\]
and $\ker(D^*) = \operatorname{Span}((I_n,\dots,I_n))$ using~\cref{eq:D1star}, which is a one-dimensional space.
Thus, the Schur matrix needs to be first projected onto its range to solve the system.
A second difficulty is that even if a sparse approximation is chosen for the Hessian~$\hat{H}$, the Schur matrix may be dense, and badly conditioned. 
However, computing the Schur matrix would be mandatory to build a preconditioner for the conjugate-gradient method for~\cref{eq:schur}. 

Regarding the iterative resolution of~\cref{eq:saddle_point}, several preconditioners have been studied and proposed in~\cite{facca24} in the context of classical optimal transport. However it is not so clear how to adapt these for the dynamical QOT distance, hence this is left for future work.

After careful comparisons between these methods the most efficient ones turned out to be the following.
To get a coarse solution, we use the Schur complement method with a block-diagonal approximation of the Hessian~\cref{eq:discrete_hessian_blockdiag},  whose inverse can therefore be computed efficiently.
When the block-diagonal approximation of the Hessian is not good enough to ensure the convergence, the true Hessian~\cref{eq:discrete_hessian} is used, and the system~\cref{eq:saddle_point} is solved either directly or with a Schur complement method. 
Both methods perform similarly.

\section{Numerical experiments} \label{sec:numerics}
In this section, we perform a series of numerical experiments and explore the performance of the algorithm. We provide the Julia source code to reproduce the figures of this article at the address \url{https://zenodo.org/records/20585024}.

\subsection{Performance}

\begin{figure}[b!]
    \centering
    \subcaptionbox{Legend}{\includegraphics[width=0.17\textwidth]{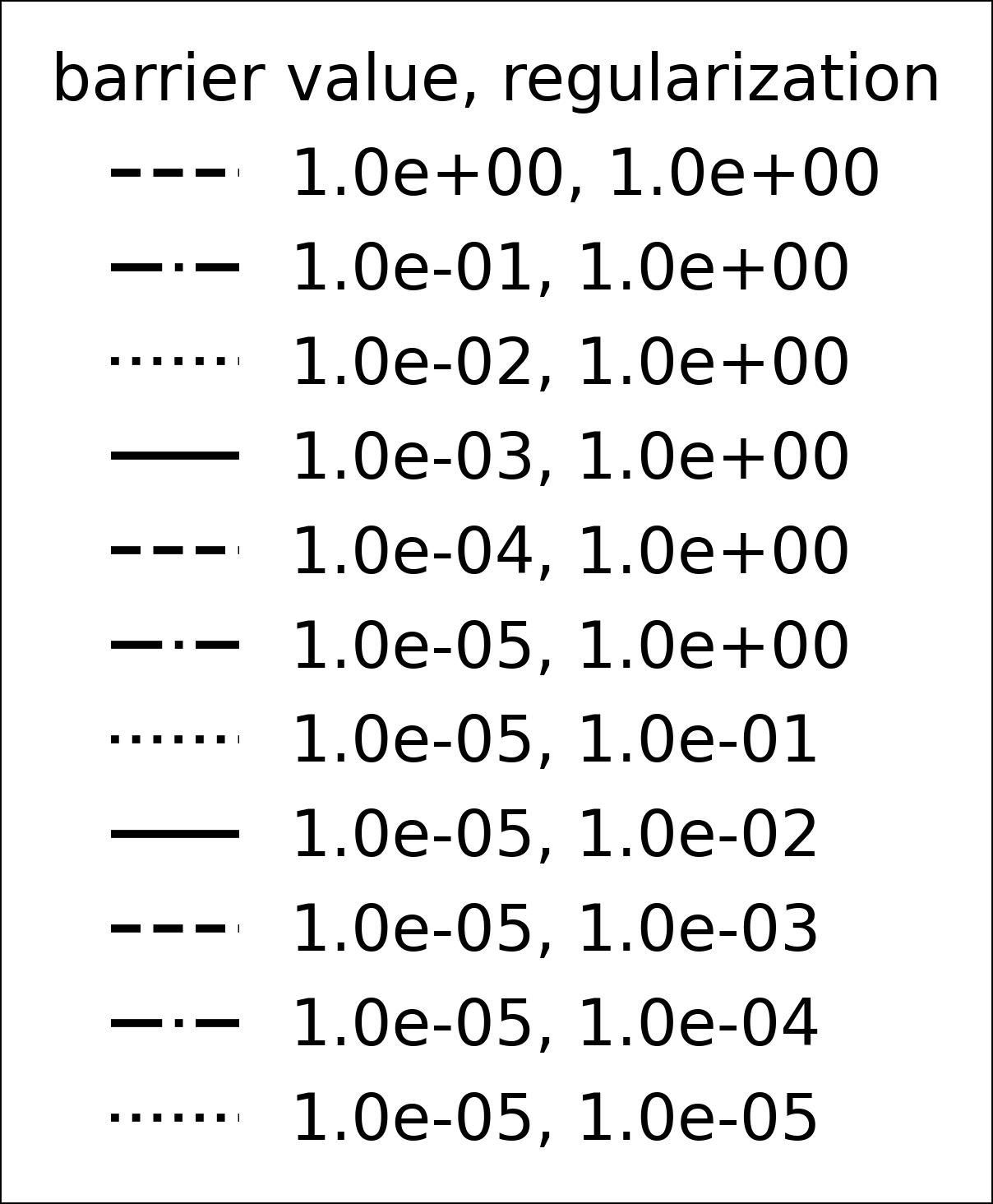}}
    \hfill
    {\includegraphics[width=0.4\textwidth]{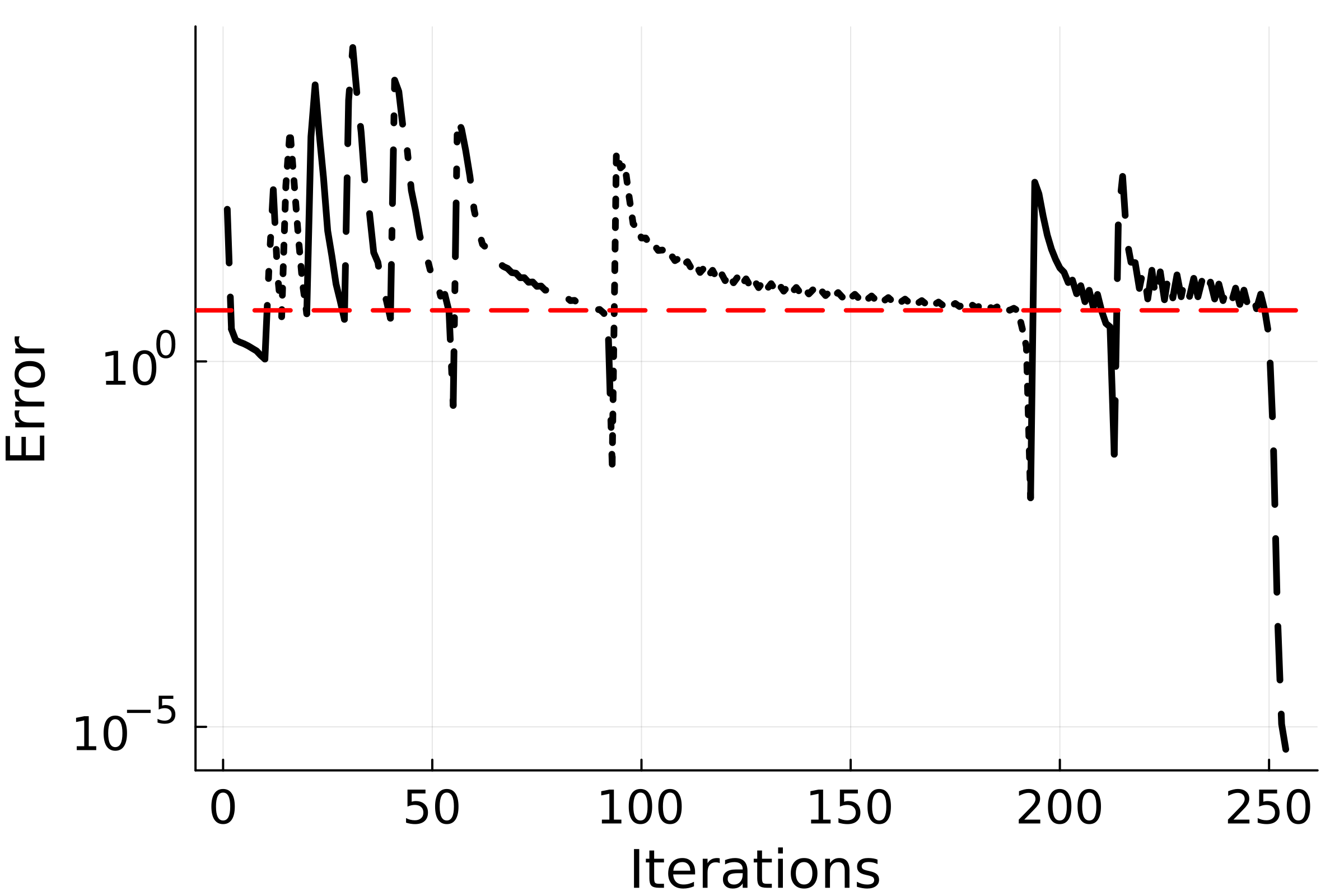}}
    {\includegraphics[width=0.4\textwidth]{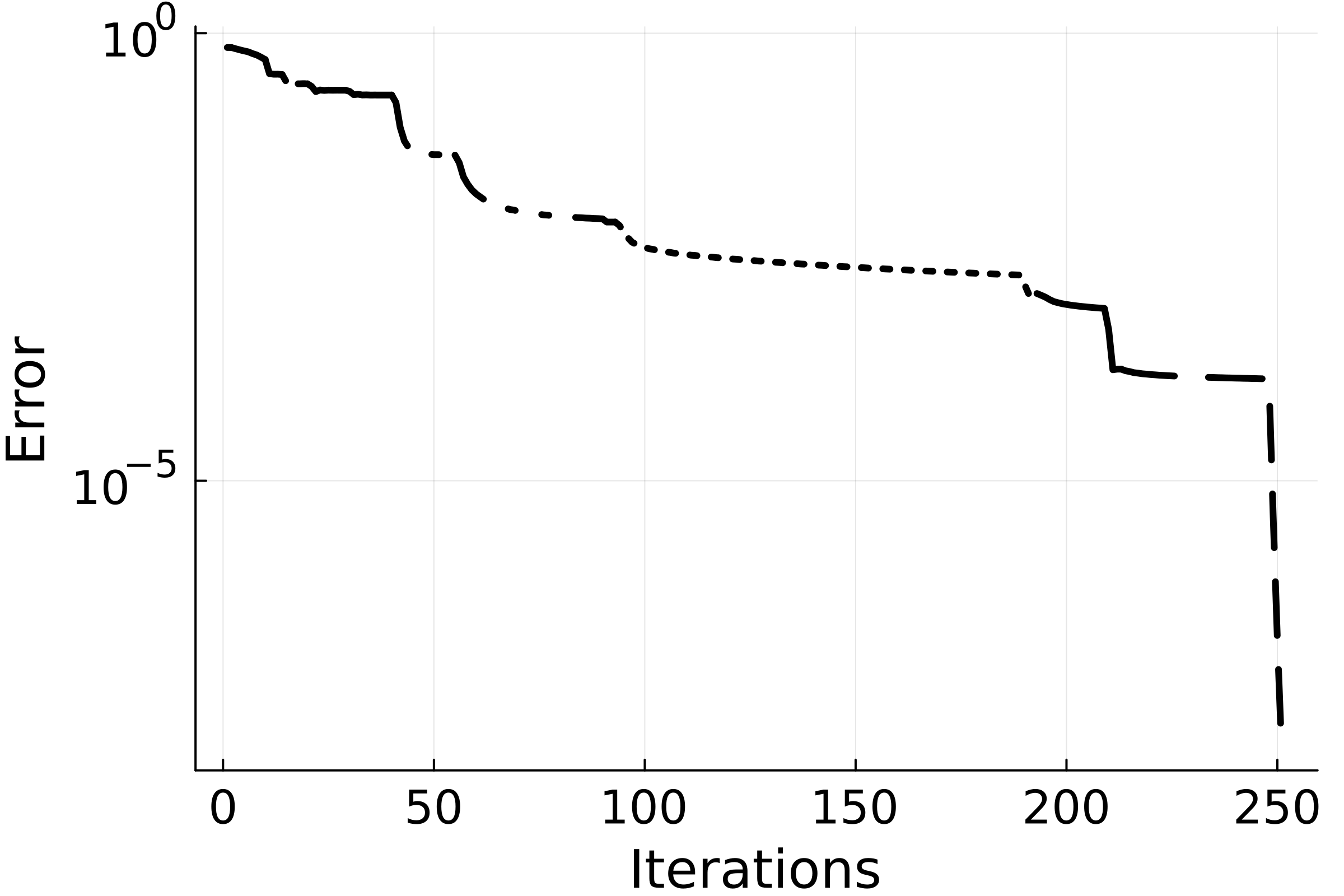}}

    \hfill
    \subcaptionbox{Residual error  $\|\nabla\mathcal{L}_\mu(\rho, \mathbf{u}, \lambda)\|_2$.}
    {\includegraphics[width=0.4\textwidth]{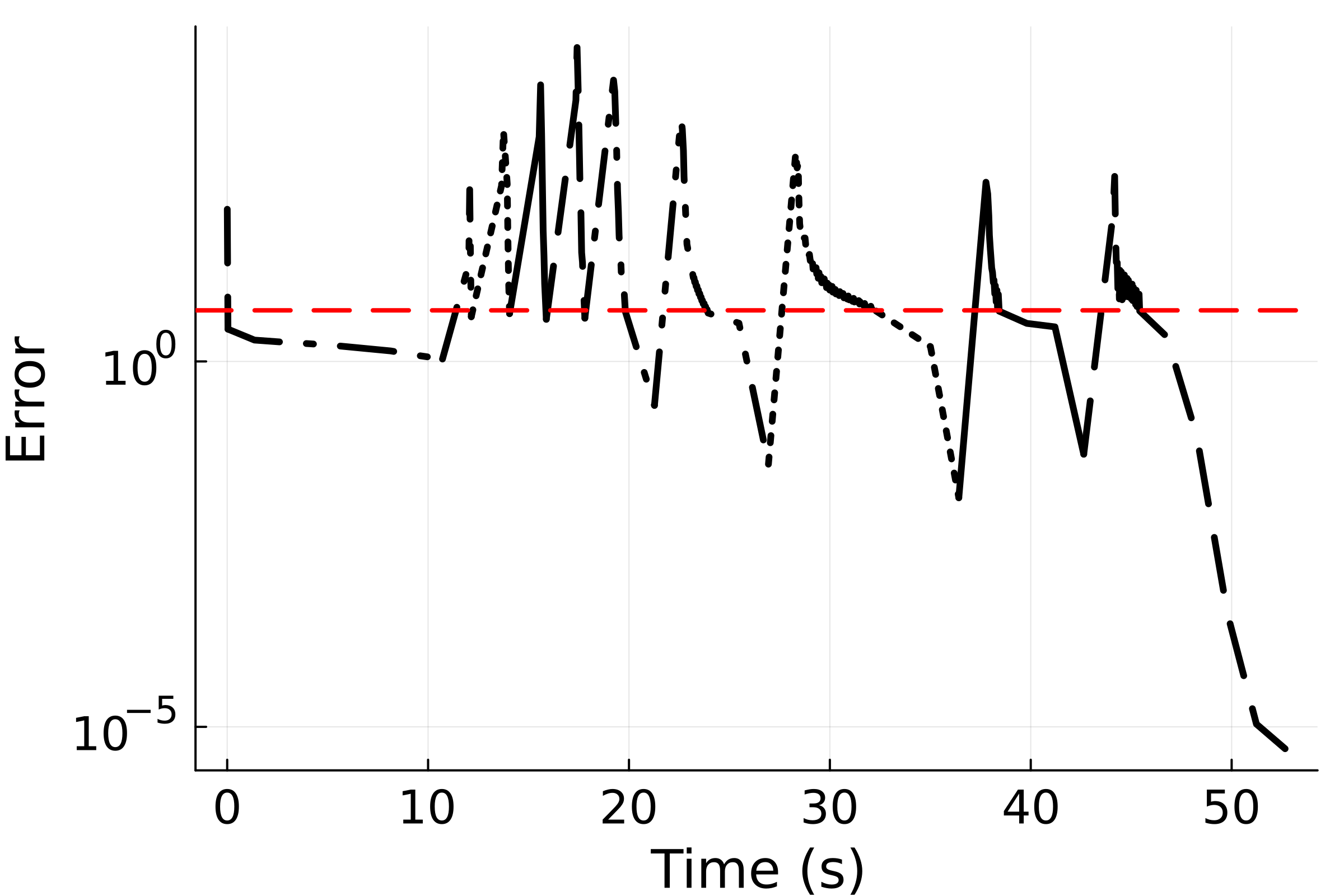}}
    \subcaptionbox{Error~\cref{eq:d_inf_2} to a reference solution.} {\includegraphics[width=0.4\textwidth]{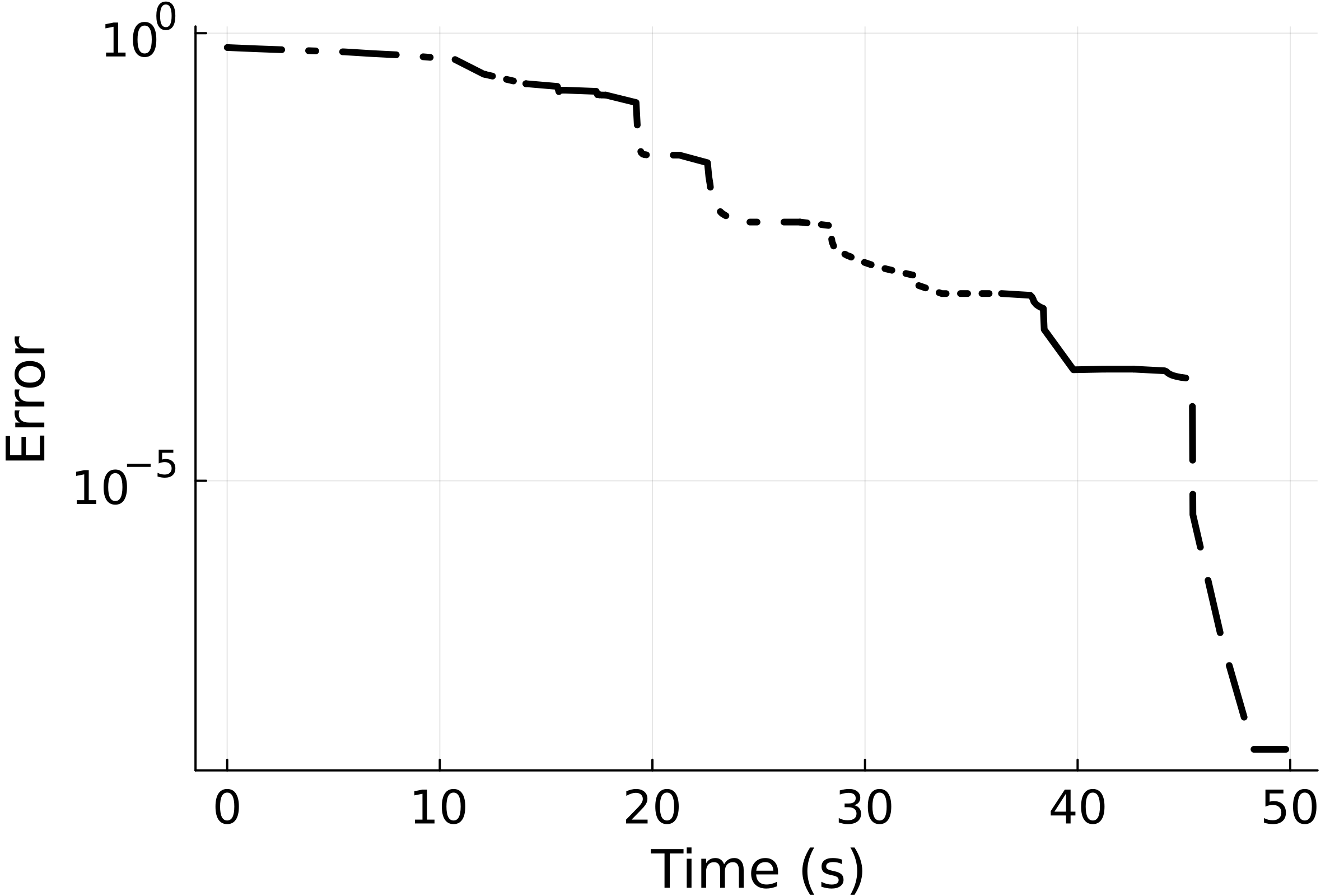}}

    \caption{Low rank computations with dimension $n=15$, $P=5$.
    The legend displays the regularization and barrier values, decreasing from $(1,1)$ to $(10^{-5},10^{-5}).$} 
    \label{fig:performance_chem_iter}
\end{figure}

We first investigate the performance of the descent~\cref{alg:descent}.
The barrier method and regularization are strongly needed for low-rank matrices, and we illustrate this fact on a test case involving two rank-two matrices $\rho_0,\rho_1$ of size fifteen arising from quantum chemistry calculations. The kernels are shown on \cref{fig:chem_geod} (respectively left and right images) and the corresponding geodesic is shown in the intermediate pictures.
In \cref{fig:performance_chem_iter} we plot on the right panel the $\ell^\infty , \ell^2$ norm to a reference solution as the number of iterations increases (top) and as a function of time (bottom). 
This error is defined by
\begin{equation}\label{eq:d_inf_2}
    \operatorname{d}_{\infty,2}(\rho,\rho^{\text{ref}}) =  \|\rho-\rho^{\text{ref}}\|_{\infty,2} 
    = \max \left\{ \| \rho_{p + \frac12}- \rho^\text{ref}_{p + \frac12} \|_2, \,  p\in \{1,\dots,P-1 \} \right\}
\end{equation}
recalling that $\|\cdot\|_2$ denotes the matrix Frobenius norm.
When $P$ and $n$ tend to infinity, this error tends to the norm
\begin{equation*}
 \|\rho-\rho^{\text{ref}}\|_{\infty,2} = \sup_{t\in[0,1]} \HS{\rho_t - \rho^{\text{ref}}_t}. 
\end{equation*}
The reference geodesic $\rho_t$ is computed such that the residual error measured as the 
$\ell^2$ norm 
of the gradient of the Lagrangian 
$ \|\nabla\mathcal{L}_{\mu}(\rho, \mathbf{u}, \lambda)\|_2  $ is lower than $7.10^{-6}.$ 

We also plot the residual error $ \|\nabla\mathcal{L}_{\mu}(\rho, \mathbf{u}, \lambda)\|_2  $ on the left of \cref{fig:performance_chem_iter}.
We observe the convergence of the error with a sharp decrease when the number of iterations is close to 250 due to a switch between approximate and exact Hessian~\cref{eq:discrete_hessian} once the residual error falls below 5.

\subsection{Convergence with respect to the size of matrices}

\begin{figure}[b!]
    \centering
 \begin{tabular}{m{1.2cm}m{10cm}}
    {\small dim.} & \\
      $n=5$ &  
  %n5
      {\includegraphics[width=.14\textwidth]{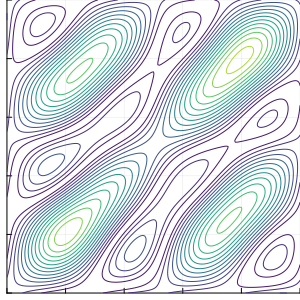}} 
      {\includegraphics[width=.14\textwidth]{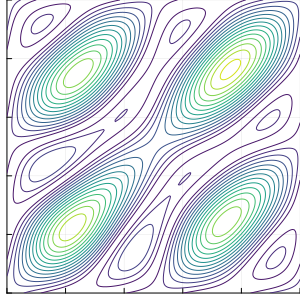}} 
      {\includegraphics[width=.14\textwidth]{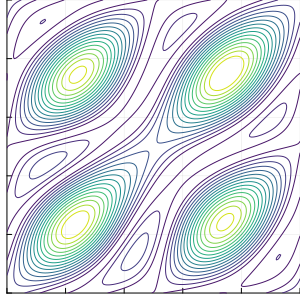}} 
      {\includegraphics[width=.14\textwidth]{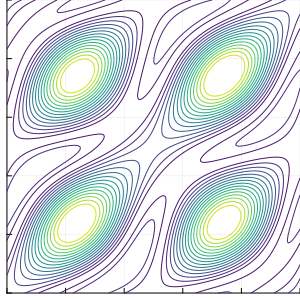}} 
      {\includegraphics[width=.14\textwidth]{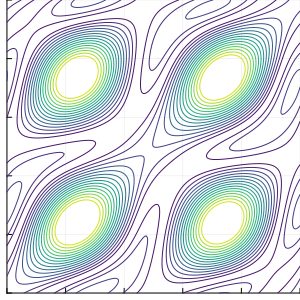}} 
  \\  
  $n=13$ &
    %n13
      {\includegraphics[width=.14\textwidth]{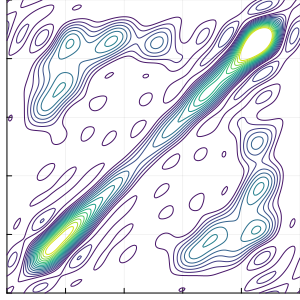}} 
      {\includegraphics[width=.14\textwidth]{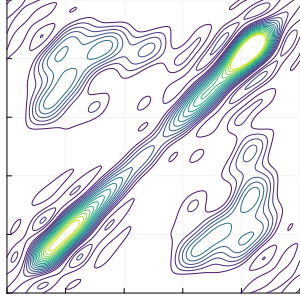}} 
      {\includegraphics[width=.14\textwidth]{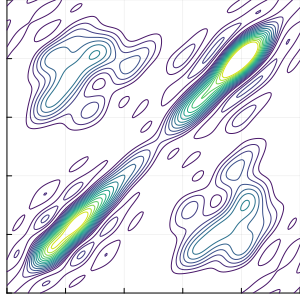}} 
      {\includegraphics[width=.14\textwidth]{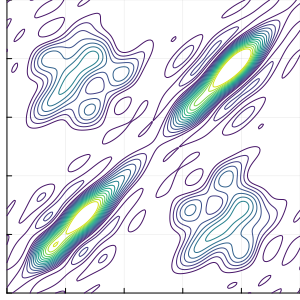}} 
      {\includegraphics[width=.14\textwidth]{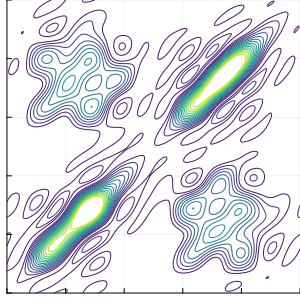}} 
  \\
  $n=21$ &
    %n21
      {\includegraphics[width=.14\textwidth]{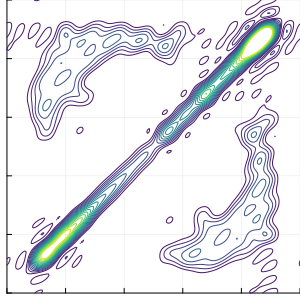}} 
      {\includegraphics[width=.14\textwidth]{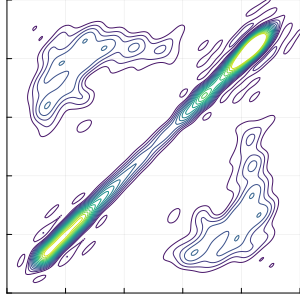}} 
      {\includegraphics[width=.14\textwidth]{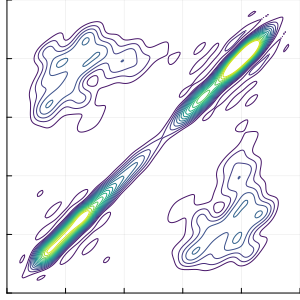}} 
      {\includegraphics[width=.14\textwidth]{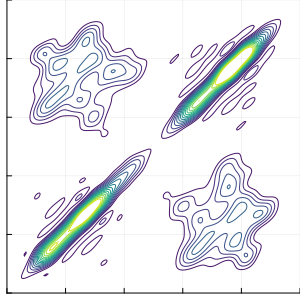}} 
      {\includegraphics[width=.14\textwidth]{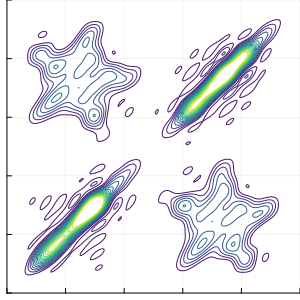}} 
  \\
  $N=31$ &
    % n31    
    \subcaptionbox{$t=0$}{\includegraphics[width=.14\textwidth]{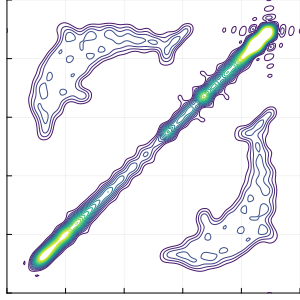}} 
    \subcaptionbox{$t=0.25$}{\includegraphics[width=.14\textwidth]{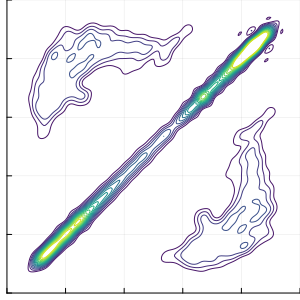}}
    \subcaptionbox{$t=0.5$}{\includegraphics[width=.14\textwidth]{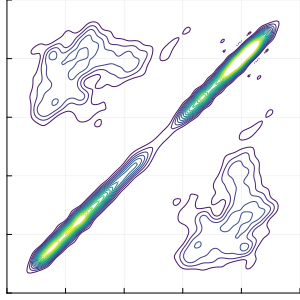}}
    \subcaptionbox{$t=0.75$}{\includegraphics[width=.14\textwidth]{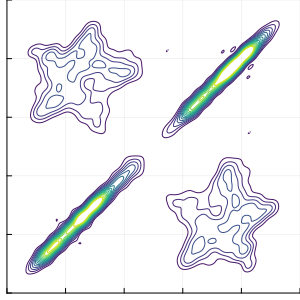}}
    \subcaptionbox{$t=1$}{\includegraphics[width=.14\textwidth]{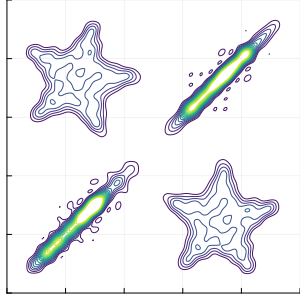}}
 \end{tabular}
\caption{Geodesics for different values of $n$. The kernels are visualized through a heatmap plot with values between 0 and 8.}
\label{fig:multiple_geods}
\vspace{-1cm}
\end{figure}

We  now numerically study the convergence of the geodesics with respect to the discretization parameter $n$.
We consider the density matrices $\rho_0^N,\rho_1^N$, with $N=31$ displayed at the bottom left and right of \cref{fig:multiple_geods} as reference density matrices. We then construct approximations thereof of size $n$ for $n=5,13,21$ by truncating the highest-frequency Fourier modes. Define $\Pi_n : \rho \in \mathcal{M}_N(\mathbb C) \mapsto  \Pi_n (\rho) =(\rho_{k,l})_{k,l=-(n-1)/2,\ldots,(n-1)/2} \in \Mn$, and 
\begin{align*}
    \rho_0^n & = \frac{\Pi_n (\rho^N_0)}{\tr(\Pi_n (\rho^N_0))} \in \Dn, \\
    \rho_1^n &= \frac{\Pi_n (\rho^N_1)}{\tr(\Pi_n (\rho^N_1))} \in \Dn.
\end{align*}

\begin{figure}[t!]
    \subcaptionbox*{}{ \includegraphics[width=0.46\linewidth]{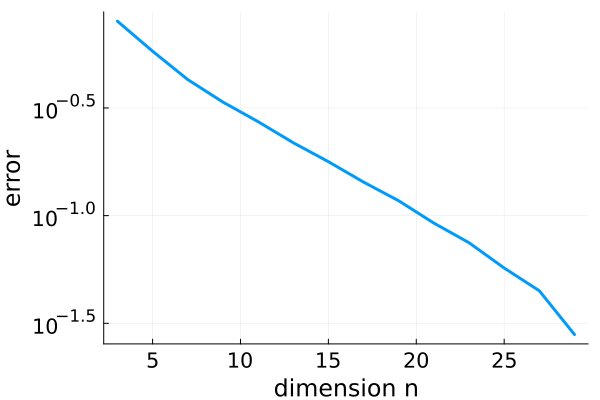}}\label{fig:err_n}
    \subcaptionbox*{}{\includegraphics[width=.46\linewidth]{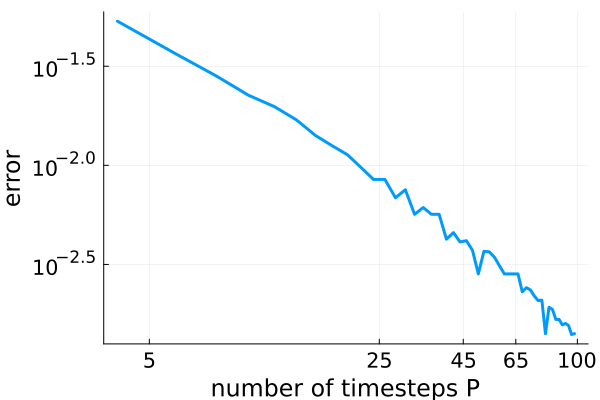}}\label{fig:err_nt}
    \caption{Errors $d_{\infty,2}$ to a reference geodesic, (left) with respect to the matrix size $n$, (right) with respect to the number of timesteps $P$.
    }
    \label{fig:cv_parameters}
    
\end{figure}

We plot the kernels over the geodesics obtained in~\cref{fig:multiple_geods} and we observe that the geodesics seem to converge when $n$ increases. To confirm this, we plot in~\Cref{fig:cv_parameters} (left) the error \cref{eq:d_inf_2} between the vector of density matrices $\rho^n$ and $\rho^N$ as a function of $n$.
Note that, for $p \in \{1,\dots,P\},$ $\rho_{p + \frac12}^n$ and $\rho_{p + \frac12}^N$ are not matrices of the same size when $n < N$. We pad the matrices $\rho_{p + \frac12}^n$ with zeros and use the quantity 
\begin{equation*}
        \left \| \rho_{p + \frac12}^N - \begin{pmatrix}
    (0) & & \\ 
    & \rho_{p + \frac12}^n& \\
    & & (0)
    \end{pmatrix} \right\|_2  
\end{equation*}
to evaluate the difference between the two, which is equivalent to computing the distance between the associated integral kernels $ \left\| \gamma_{p + \frac12}^n - \gamma_{p + \frac12}^N \right\|_2$.
The observed error decreases at the same rate as the error 
\[
    \max\left(\| \rho_0^n - \rho_0^N\|_2, \, \| \rho_1^n - \rho_1^N\|_2\right)
\]
computed with padded zeros. This observation supports the idea that the geodesics converge as $n \to +\infty$. 

\subsection{Convergence as the number of timesteps grows}

We then study the convergence of the dynamical QOT distance with respect to the number of timesteps $P$, in order to find the number of steps suited for the discretization of the problem at hand.
To explore this numerically, we compute a reference solution with $P_{\rm ref}=100$ and compare it to solutions computed with varying time steps $P$. All density matrices are of size $n=7$ and $\rho_0^n$ and $\rho_1^n$ correspond to the ones on \Cref{fig:multiple_geods}.
Since the solutions $(\rho_p)$ do not have the same dimension (vectors of size $P$ vs vectors of size $P_{\rm ref}$), 
we extend the vectors to piecewise constant vectors to ensure that every discrete curve has the same number of time samples. 
More precisely, if the length of $\rho_1$ is ${P}_1$ and the length of $\rho_2$ is ${P}_2$, with ${P}_1 < {P}_2,$ we actually compare $\rho_2$ with the vector $\widetilde{\rho}_1$ of length ${P}_2$, whose $p^\text{th}$ element is defined by
\begin{equation}
\label{eq:rhoPapprox}
    (\widetilde{\rho}_1)_p = ({\rho_1})_{\tilde{p}} \text{, with } \quad  \tilde{p} = \left\lfloor (p-1) \frac{{P}_1}{{P}_2}\right\rfloor +1.
\end{equation}
We see in~\Cref{fig:cv_parameters} (right) that the error decreases when $P$ increases, and we also observe that this error is comparable to the projection error of the geodesic of length $P_{\rm ref}$ to a geodesic of length $P$ using~\cref{eq:rhoPapprox}. 

\subsection{Comparison with classical optimal transport}\label{subsec:comparisonOT}

\begin{figure}[b!]
    \centering
    % our geodesic
    \includegraphics[width=0.19\textwidth]{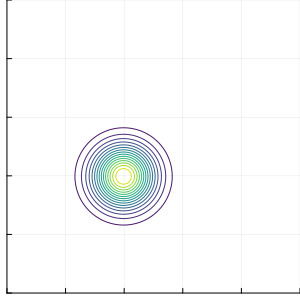}
    \includegraphics[width=0.19\textwidth]{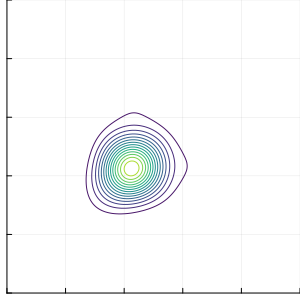}
    \includegraphics[width=0.19\textwidth]{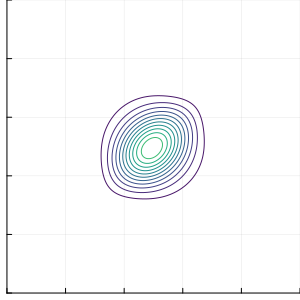}
    \includegraphics[width=0.19\textwidth]{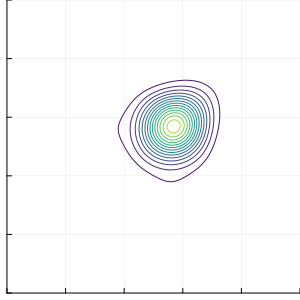}
    \includegraphics[width=0.19\textwidth]{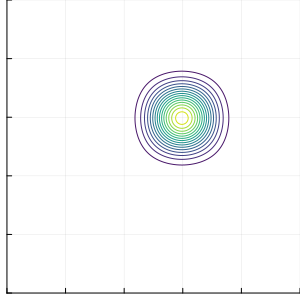}
    
    %ot 
    \subcaptionbox{$t=0$}{\includegraphics[width=0.19\textwidth]{figures/fig_10_-_OT_vs_QOT_2D__2_/OT_kernels_2_t0.0.png}}
    \subcaptionbox{$t=0.25$}{\includegraphics[width=0.19\textwidth]{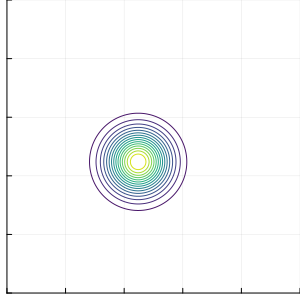}}
    \subcaptionbox{$t=0.5$}{\includegraphics[width=0.19\textwidth]{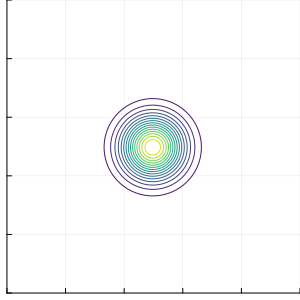}}
    \subcaptionbox{$t=0.75$}{\includegraphics[width=0.19\textwidth]{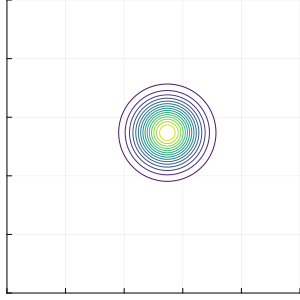}}
    \subcaptionbox{$t=1.0$}{\includegraphics[width=0.19\textwidth]{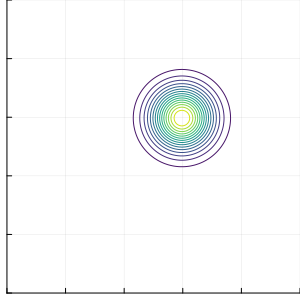}}
    \caption{Geodesic between translated kernels (heatmaps with values between 0 and 8). Top:  dynamical QOT geodesic. Bottom: classical 2D optimal transport geodesic.}
    \label{fig:translation2}
\end{figure}

\begin{figure}[b!]
    \centering
    \subcaptionbox{Densities along a dynamical QOT geodesic for Gaussian kernels.}{\includegraphics[width=0.49\linewidth]{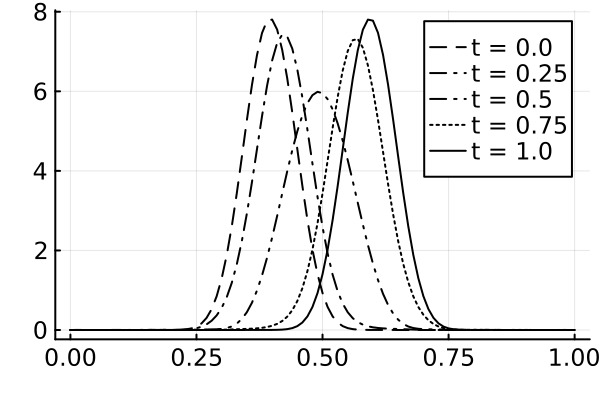}}
    \subcaptionbox{A classical 1D OT geodesic for Gaussian functions.}{\includegraphics[width=0.49\linewidth]{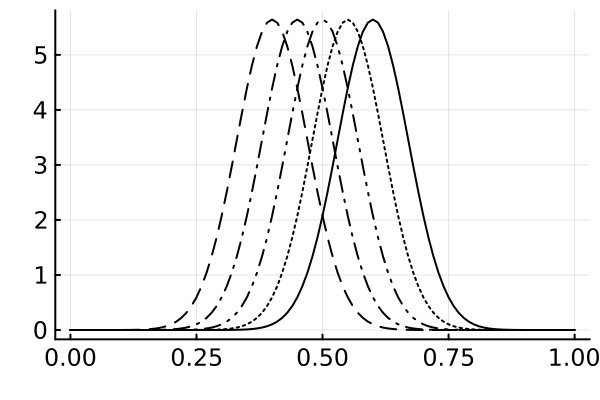}}
    \caption{Comparison between the diagonal in the first line of~\cref{fig:translation2} and a classical optimal transport geodesic.}
    \label{fig:densities_translation2}
\end{figure}

Although the dynamical formulation of QOT is \emph{a priori} not directly linked with classical OT, it is nevertheless possible to compare the two formulations in at least two ways. 
First it is possible to compare the diagonals of kernels for dynamical QOT geodesics which are probability measures with corresponding classical OT geodesics.
Second we can compare dynamical QOT geodesics and classical OT geodesics between  density matrices with everywhere positive kernels of equal mass.
We provide these comparisons with Gaussian densities.

\begin{figure}[b!]
    \centering
    % our geodesic
    \includegraphics[width=0.19\textwidth]{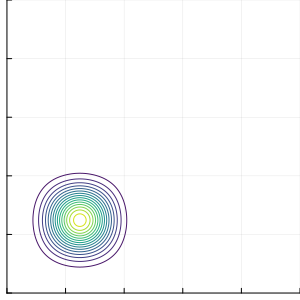}
    \includegraphics[width=0.19\textwidth]{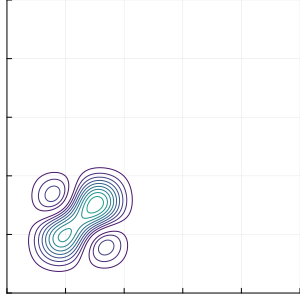}
    \includegraphics[width=0.19\textwidth]{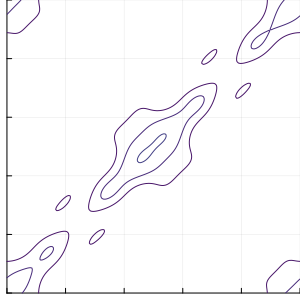}
    \includegraphics[width=0.19\textwidth]{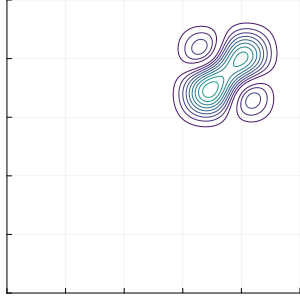}
    \includegraphics[width=0.19\textwidth]{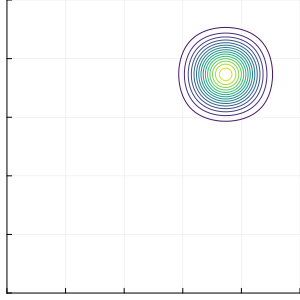}
    
    \subcaptionbox{$t=0$}{\includegraphics[width=0.19\textwidth]{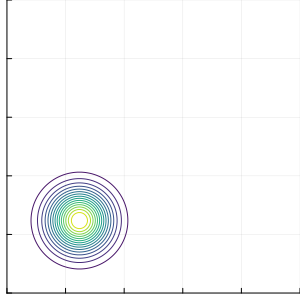}}
    \subcaptionbox{$t=0.25$}{\includegraphics[width=0.19\textwidth]{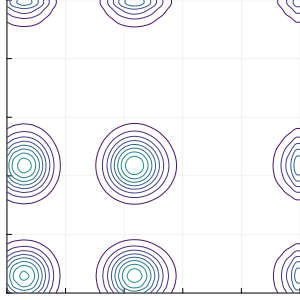}}
    \subcaptionbox{$t=0.5$}{\includegraphics[width=0.19\textwidth]{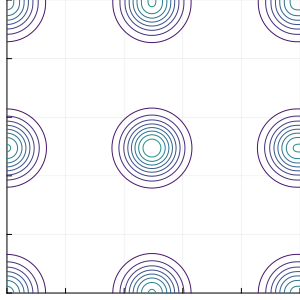}}
    \subcaptionbox{$t=0.75$}{\includegraphics[width=0.19\textwidth]{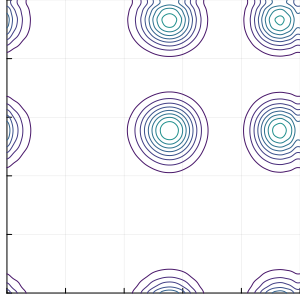}}
    \subcaptionbox{$t=1.0$}{\includegraphics[width=0.19\textwidth]{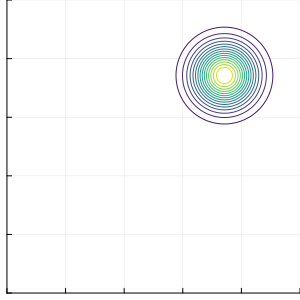}}
    \caption{Geodesic between Gaussian kernels (heatmaps with values between 0 and 8). Top: dynamical QOT geodesic. Bottom: classical optimal 2D transport geodesic.}
    \label{fig:translation}
\end{figure}

\begin{figure}[b!]
    \centering
    \subcaptionbox{Densities along a dynamical QOT geodesic for Gaussian kernels.}{\includegraphics[width=0.49\linewidth]{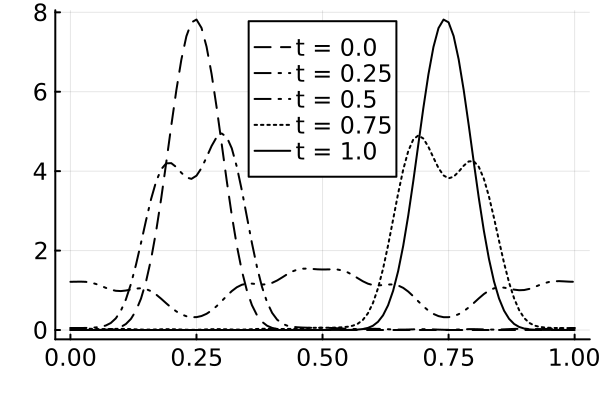}}
    \subcaptionbox{A classical 1D OT geodesic for Gaussian functions.}{\includegraphics[width=0.49\linewidth]{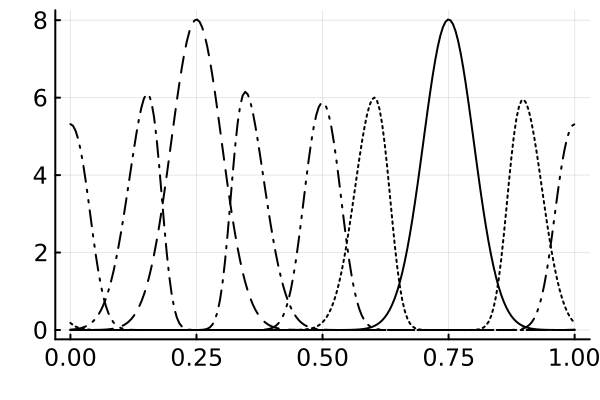}}
    \caption{Comparison between the diagonal in the first line of~\cref{fig:translation} and a classical optimal transport geodesic.}
    \label{fig:densities_translation}
\end{figure}

We consider $\gamma_0$, $\gamma_1$ two Gaussian functions with means
$\mu_0 = (2/5, 2/5)$ and $\mu_1 = (3/5,3/5)$ and covariance matrices $\Sigma_0 = \Sigma_1 = 5. 10^{-3} I_2$.
We approximate the kernels $\gamma_0$ and $\gamma_1$ with a Fourier discretization of size 13, see~\cref{eq:gam_rho}.
We compute the QOT geodesic between the density matrices associated to the discretized kernels using a lowest regularization of $\varepsilon = 10^{-5}$ and a lowest barrier of $\mu=10^{-5}$.
We plot the geodesic on~\cref{fig:translation2} (top), and we include the classical OT geodesic between $\gamma_0$ and $\gamma_1$ for comparison (bottom), which is analytically known.
We observe in the two cases that the mass translates. 
However, in classical OT, the mass stays concentrated while in QOT computations, the Gaussian densities seem to flatten, as in regularized classical OT computations.
This can be clearly observed on~\cref{fig:densities_translation2}, where we plot the diagonals $\gamma_t(x,x)$ of the intermediate QOT density matrices,  as well as the classical OT geodesic between the 1D Gaussians $\gamma_0(x,x)$ and $\gamma_1(x,x)$ which are just translated Gaussians.

We then consider two different Gaussian kernels $\gamma_0$, $\gamma_1$ with means
$\mu_0 \hspace{-1pt} = (1/4, 1/4)$, $\mu_1 = (3/4,3/4)$ and covariance matrices $\Sigma_0 = \Sigma_1 = 5. 10^{-3} I$.
We compute the QOT geodesic between the density matrices associated to the discretized kernels with $n=13$.
We plot the geodesic on~\cref{fig:translation} (top), and we include the classical OT geodesic between $\gamma_0$ and $\gamma_1$ for comparison (bottom).
The classical OT geodesic is computed with the Python Optimal Transport (POT) library~\cite{flamary21}, using  the \verb|bregman.barycenter()| function~\cite{benamou15}.
We observe in the two cases that the mass both splits due to periodic boundary conditions and translates. 
In the QOT geodesic, we observe a clear regularizing behavior, which is not present for the classical OT case. As before, this can be easily observed when considering the diagonals on~\cref{fig:densities_translation2}.

\subsection{Rank/eigenvalue evolution along a geodesic}\label{subsec:rank}
\begin{figure}[htb!]
    \centering
    \subcaptionbox{Eigenvalues evolution in the example of~\cref{fig:translation}.}{\includegraphics[width=0.49\textwidth]{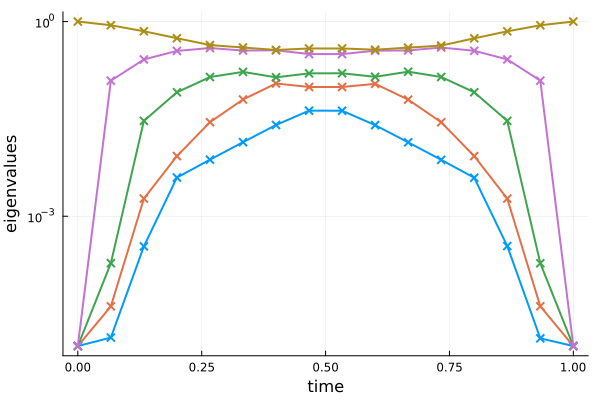}}
    \subcaptionbox{Eigenvalues evolution in the first line of the example in~\cref{fig:chem_geod}.}{\includegraphics[width=0.49\textwidth]{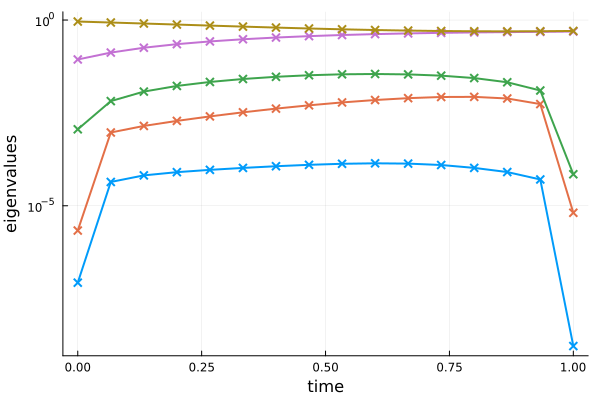}}
    \caption{Smallest eigenvalues along dynamical QOT geodesics, showing that intermediate points are less singular than the endpoints ($n=5$).}
    \label{fig:eigen_variation}
\end{figure}

The dynamical QOT distance~\cite{carlen20} is primarily defined on positive \emph{definite} matrices. To observe if the intermediate points in the geodesics are less singular than the endpoints, we plot on~\cref{fig:eigen_variation} the lowest eigenvalues of the density matrices along two different geodesics appearing respectively in \Cref{fig:chem_geod} and \Cref{fig:translation}. They are several orders of magnitude larger in intermediate points of the geodesic than in endpoints, which 
seems to indicate that the low-rank of the endpoints is not preserved by the geodesics.

\section{Quantum chemistry experiments and choice of derivations} \label{sec:chem}

\begin{figure}[htb!]
    \centering
    % QOT
    \includegraphics[width=0.19\textwidth]{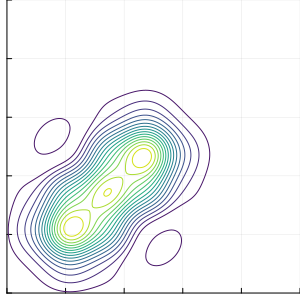}
    \includegraphics[width=0.19\textwidth]{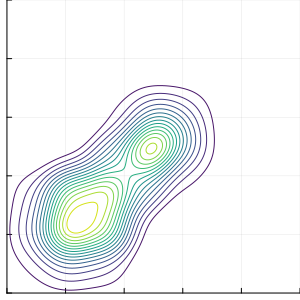}
    \includegraphics[width=0.19\textwidth]{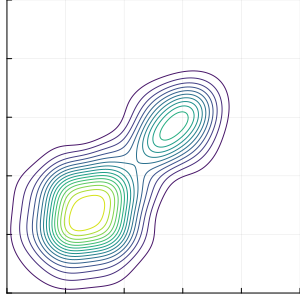}
    \includegraphics[width=0.19\textwidth]{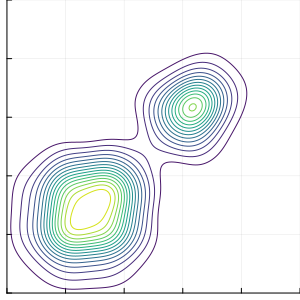}
    \includegraphics[width=0.19\textwidth]{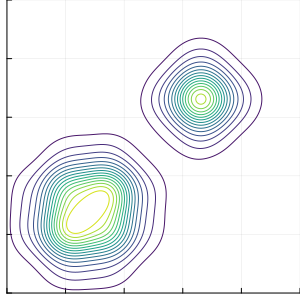}
    
    % FE 
    \subcaptionbox{$t=0$}{\includegraphics[width=0.19\textwidth]{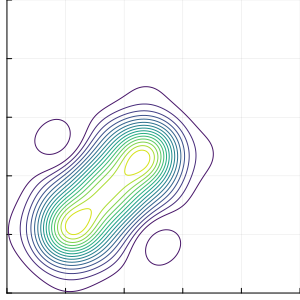}}
    \subcaptionbox{$t=0.25$}{\includegraphics[width=0.19\textwidth]{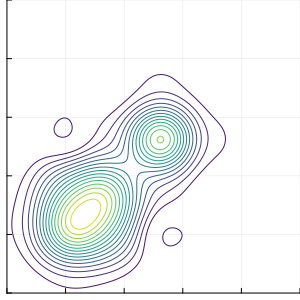}}
    \subcaptionbox{$t=0.5$}{\includegraphics[width=0.19\textwidth]{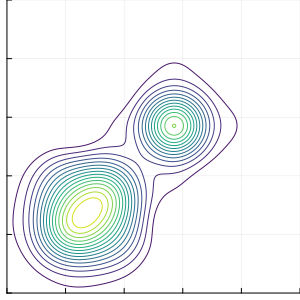}}
    \subcaptionbox{$t=0.75$}{\includegraphics[width=0.19\textwidth]{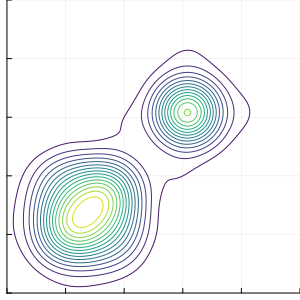}}
    \subcaptionbox{$t=1.0$}{\includegraphics[width=0.19\textwidth]{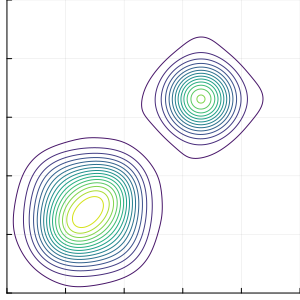}}
    \caption{Comparison between the best-fitting QOT geodesic (Top) and the reference computation (Bottom), a three-electron simulation.}
    \label{fig:translation_three_electrons}
\end{figure}

\begin{figure}[htb!]
    \centering
    \subcaptionbox{Best-fitting QOT geodesic densities.}{\includegraphics[width=0.49\linewidth]{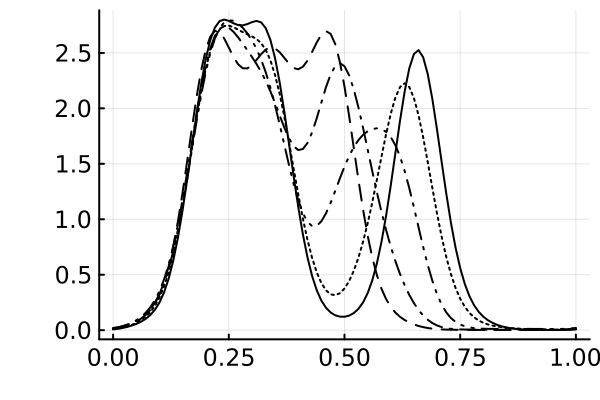}}
    \subcaptionbox{Quantum chemistry densities.}{\includegraphics[width=0.49\linewidth]{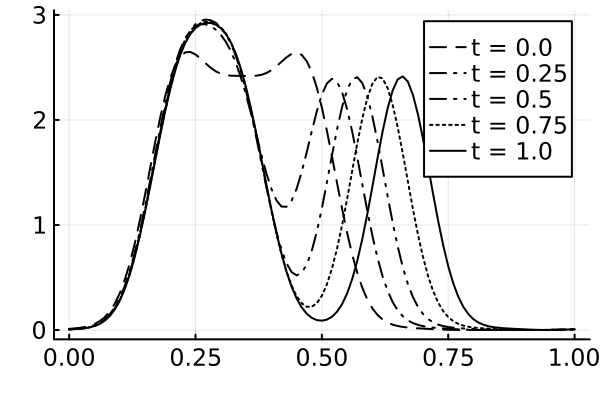}}
    \caption{Comparison between densities in the three-electron experiment.}
    \label{fig:densities_translation_three_electrons}
\end{figure}

One of the main motivations of this article is to 
compare the QOT geodesics with natural curves of density matrices coming from quantum chemistry.
More precisely, we consider, as in~\Cref{subsec:motivation}, different density matrices $T_{\gamma_\theta}$ that come from the solution of the electronic Schrödinger equation for varying nuclei positions.
We compute these density matrices using the Julia package \texttt{SchrodingerFE.jl}~\cite{quan23}. 
For two-electron systems we use two nuclei each of charge 1 with distances $0.15, 0.213, 0.275, 0.338, 0.4$.
For three electron systems  we consider for the nuclei two positive charge 1 at positions $0.25$ and $0.375$, and another  positive charge 1  at positions $0.5, 0.545, 0.588, 0.631,$ and $0.675$.
As shown on \cref{fig:chem_geod} and \cref{fig:densities_translation_three_electrons}
the chemistry and QOT results look very similar.
Note that the derivation operators are chosen differently in the two and three electrons cases: 
for two electrons we use $L_1(0.83)$ and $L_2(-\infty)$ while for three electrons we use $L_1(2)$ and $L_2(-\infty)$. 
These geodesics correspond in fact to optimal parameter choices as we detail below.

\begin{figure}[htb!]
    \centering
    \subcaptionbox{Two-electron case.
    }{\includegraphics[width=.45\linewidth]{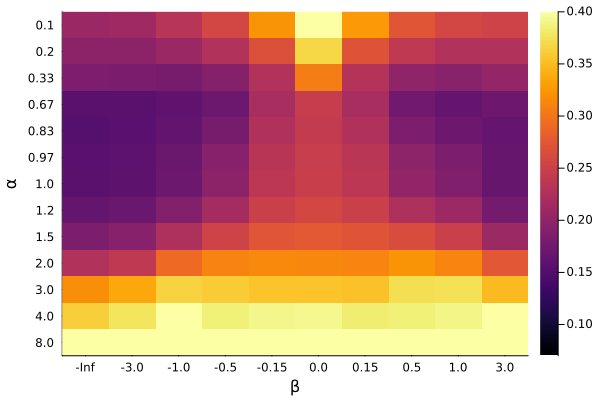}}    
    \subcaptionbox{Three-electron case.
    }{\includegraphics[width=.45\linewidth]{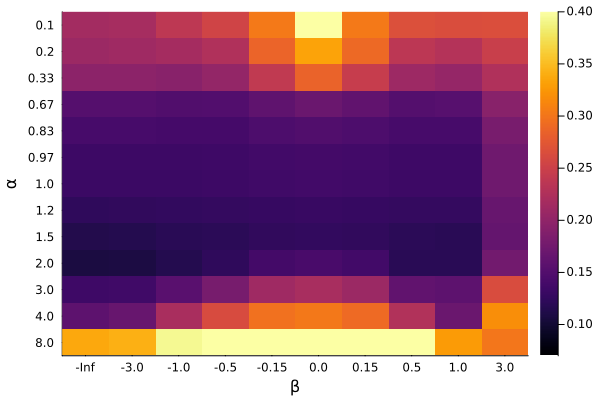}}
    \caption{Heatmaps of errors $d_{\infty,2}$ of a set of QOT geodesics with respect to a reference computation from quantum chemistry, with varying derivation parameters $\alpha$ and $\beta$.}
    \label{fig:error_heatmap}
\end{figure}

Since the dynamical QOT distance depends on derivation choices, we test multiple QOT distances for the derivation set defined in~\crefrange{eq:derivations1}{eq:derivations2} and compare them with quantum chemistry results. 
We therefore test a large parameter set $(\alpha,\beta)$ defining $L_1(\alpha)$ and $L_2(\beta)$ and compute the error $d_{\infty,2}$ (see~\cref{eq:d_inf_2}) between the QOT geodesics and quantum chemistry reference curves displayed on~\cref{fig:error_heatmap}.
We therefore observe that the optimal parameters are 
 $(\alpha=0.83, \beta=-\infty)$ for the two-electron case, and $(\alpha=2, \beta=-\infty)$ for the three-electron case.
At this stage we do not have an explanation for these optimal parameters but understanding this phenomenon seems a promising research direction; we plan to explore this in the future.
Note that while the geodesics on \Cref{fig:translation_three_electrons} seem very close the diagonals plotted on \Cref{fig:densities_translation_three_electrons} are clearly different. Indeed a larger regularizing effect is still observed on the dynamical QOT geodesics compared to the quantum chemistry electronic densities. 

\section{Conclusion and future work}
\label{sec:concl}
In this work, we developed an interior-point, regularized method based on the  SQP algorithm introduced in~\cite{chen17}. We showed that we could carry out computations for singular positive semidefinite / low-rank density matrices, and used this method to investigate whether the dynamical QOT distance introduced in~\cite{carlen20} is appropriate to approximate solutions in quantum chemistry. It seems that QOT geodesics can provide a good fit to curves of density matrices arising in quantum chemistry, provided that the derivations and functional calculus are chosen appropriately. 
This work opens the way to many interesting research directions such as the practical computation of generic QOT Wasserstein barycenters based on this dynamical formulation, and the construction of low-rank dynamics for density matrices, which would alleviate the memory requirements of our algorithm and be easily usable for low-rank quantum chemistry calculations.

\section*{Acknowledgements}
This work has been supported by the EIPHI Graduate school (contract ANR-17-EURE-0002) and by the Région 
Bourgogne-Franche-Comté.
This work has received funding from the ANR through the project NUMERIQ (grant number ANR-24-CE46-2255). Etienne Obermeyer thanks the Program QuanTEdu-France
n°ANR-22-CMAS-0001 France 2030. We acknowledge the financial support of European Research Council (ERC) under the European Union’s
Horizon 2020 Research and Innovation Programme – Grant Agreement n°101077204 HighLEAP. 

\medskip

We thank Yann Brenier, as well as Andrea Natale and Gabriele Todeschi for stimulating discussions.

\FloatBarrier

\appendix
\section{Hessian of the discrete objective}
\label{app:Hessian}
Let us derive the Hessian of the discrete function
\begin{align*} 
 f_\mu(\rho, \mathbf{u}) =    \left\langle \mathbf{u}^* \circ \mathbf{u},\; A(\rho^{-1}) + a \right\rangle + \mu \operatorname{Barr}(\rho),
\end{align*}
 for $(\rho, \mathbf u) \in \Dn^{P-1} \times (\Mn^{J})^P$, from which the Hessian of the corresponding function without the barrier can easily be obtained.
 Let  $g :\rho \mapsto \rho^{-1}$, and $X \in \Mn$. For an invertible matrix $\rho$, we have
$\nabla_\rho \,g(\rho)(X) = - \rho^{-1} X \rho^{-1}$, so that  
\begin{align*}
\nabla_\rho {f}_\mu(\rho,\mathbf{u}) &= -\rho^{-1} \circ A^*(\mathbf{u}^* \circ \mathbf{u}) \circ \rho^{-1} - \mu\rho^{-1}, \\
\nabla_{\mathbf{u}} f_\mu(\rho,\mathbf{u}) &=   2  \mathbf{u} \circ (A(\rho^{-1}) + a).
\end{align*}
For all $p \in \{1,\dots, P-1\}$, we define the operator $H^{\rho \rho}_{\mu,\rho_{p+\frac12},\mathbf u} : \Herm \to \Herm$ by
\begin{multline*}    
H^{\rho \rho}_{\mu,\rho_{p+\frac12},\mathbf u} (X) = \\ \rho_{p+\frac12}^{-1} A^*(\mathbf{u}^* \circ \mathbf{u}) \rho_{p+\frac12}^{-1} X \rho_{p+\frac12}^{-1} +  \rho_{p+\frac12}^{-1} X \rho_{p+\frac12}^{-1}  A^*(\mathbf{u}^* \circ \mathbf{u}) \rho_{p+\frac12}^{-1} + \mu \rho_{p+\frac12}^{-1} X \rho_{p+\frac12}^{-1}.
\end{multline*}
The operator $\nabla_\rho \nabla_\rho f_\mu(\rho, \mathbf{u}) : (\Herm)^{P-1} \to (\Herm)^{P-1}$ can be expressed as
\[
    \nabla_\rho \nabla_\rho f_\mu(\rho, \mathbf{u})(X) = \left( H^{\rho \rho}_{\mu,\rho_{p+\frac12},\mathbf u} (X_p) \right)_{p=1}^{P-1}.
\]
The operator $\nabla_\mathbf{u} \nabla_\mathbf{u} f_\mu(\rho, \mathbf{u}) : (\Mn^J)^P \to (\Mn^J)^P$ writes
\begin{align*}
    \nabla_\mathbf{u} \nabla_\mathbf{u} f_\mu(\rho, \mathbf{u})(\mathbf{X}) 
    &= 2\mathbf{X} \circ (A(\rho^{-1}) +a) \\
    &= \left( 2\underbrace{\mathbf{X}_p}_{\in \Mn^J} \underbrace{(A(\rho^{-1})+a)_p}_{\in \Herm }\right)_{p=1}^P \\
    &= \big( \left( 2(X_p)_j {(A(\rho^{-1})+a)_p} \right)_{j=1}^J\big)_{p=1}^P.
\end{align*}
Finally, one can compute $\nabla_\mathbf{\rho} \nabla_\mathbf{u} f_\mu(\rho, \mathbf{u}) : \Herm^{P-1} \to (\Mn^J)^P$, with
\[
     \nabla_\mathbf{\rho} \nabla_\mathbf{u} f_\mu(\rho, \mathbf{u})(X) = -2\mathbf{u} \circ A(\rho^{-1} X\rho^{-1}).
\]

\bibliographystyle{siam}
\bibliography{dynamicalQOT}

\end{document}